\documentclass[11pt]{article}

\usepackage{amsmath}
\usepackage{amscd}
\usepackage{latexsym}
\usepackage{amssymb}
\usepackage[all]{xy}

\newfont{\gothique}{eufm10 scaled 1100}  
\newcommand{\goth}[1]{\mbox{\gothique{#1}}}

\newcommand{\PP}{{\mathbb{P}}}
\newcommand{\CC}{{\bf{C}}}
\newcommand{\ZD}{{\cal{Z}}}
\newcommand{\JJ}{{\cal{J}}}
\newcommand{\JG}{{\bf{J}}_{\Gamma}}

\newcommand{\OO}{{\cal O}}
\newcommand{\GS}{{ \cal G}}
\newcommand{\GT}{{\tilde{ \cal G}}}

\newcommand{\EZ}{Ext_{Z}^{1}}
\newcommand{\SE}{{\cal E}}
\newcommand{\SEE}{{\SE},[e]}

\newcommand{\ENDO}{{\cal E}{\it{nd}}}

\newcommand{\JA}{{\bf{J}}( X;L,d)} 
\newcommand{\JAA}{{\bf{J}}}

\newcommand{\XD}{X^{[d]}}
\newcommand{\HT}{{\bf \tilde{H}}}
\newcommand{\HH}{{\bf{H}}}
\newcommand{\ID}{{\cal{I}}}

\newcommand{\PIPO}{\pi^{\ast}\big( p_{2 \ast} \OO_{\ZD}\big)}

\newcommand{\FIT}{{\bf{\tilde{ F}}}}
\newcommand{\FI}{{\bf F}}
\newcommand{\GA}{\Gamma}
\newcommand{\GAC}{\Gamma^{(0)}_{conf}}

\newcommand{\GG}{{\Gamma}^r_d}

\newcommand{\JAB}{{\bf \breve{J}}} 
\newcommand{\HOM}{{\cal H}\it{om}}

\newcommand{\QB}{{\bf q}}
\newcommand{\QT}{{\bf {\tilde{q}}}}

\newcommand{\HO}{H^0(\OO_}

\newcommand{\ZA}{([Z],[\alpha])}
\newcommand{\FT}{\tilde{\cal F}}
\newcommand{\FF}{{\cal F}}
\newcommand{\CS}{C^r (L,d)}
\newcommand{\TE}{{\bf \Theta}} 
\newcommand{\RI}{\cite{[R1]}}
\newcommand{\CO}{Conf_d (X)}
\newcommand{\CSA}{C^r_{adm} (L,d)}
\newcommand{\LG}{l_{\GA}}
\newcommand{\LAG}{\GS_{\mbox{\unboldmath$\GA$}} }
\newcommand{\LAGT}{\GT_{\mbox{\unboldmath$\GA$}} }
\newcommand{\CG}{{\cal C}_{\mbox{\unboldmath$\GA$}} }
\newcommand{\BM}{\boldmath}

\newcommand{\BEN}{\begin{equation}}
\newcommand{\EEN}{\end{equation}}
\newcommand{\GAB}{\breve{\GA}}
\newcommand{\JABG}{\JAB_{\GA}}
\newcommand{\FL}{{\cal{FL}}}
\newcommand{\FLA}{\mbox{\BM${\FL}$}_{h_{\GA}}}

\newtheorem{thm}{Theorem}[section]
\newtheorem{lem}[thm]{Lemma}
\newtheorem{pro}[thm]{Proposition}
\newtheorem{cor}[thm]{Corollary}
\newtheorem{rem}[thm]{Remark}

\newtheorem{defi}[thm]{Definition}  

\newtheorem{pro-defi}[thm]{Proposition-Definition}
\newlength{\myskip}
\makeatletter
\renewcommand{\@seccntformat}[1]{\S \/ {\csname the#1\endcsname}\hspace{0.5em}}
\makeatother
\title{\bf NONABELIAN JACOBIAN OF SMOOTH PROJECTIVE SURFACES - A SURVEY}
\author{Igor Reider}

\large
\begin{document}
\bibliographystyle{amsplain}
\maketitle
\numberwithin{equation}{section}
\begin{abstract}
The nonabelian Jacobian $\JA$ of a smooth projective surface $X$ is inspired
 by the classical theory of Jacobian of curves.
It is built as a natural scheme interpolating between the Hilbert scheme $\XD$ of subschemes of length
$d$ of $X$ and the stack
${\bf M}_X (2,L,d)$ of torsion free sheaves of rank $2$ on $X$ having the determinant $\OO_X (L)$ and the second Chern class (= number) $d$. It relates to such influential ideas as variations of Hodge structures, period maps, nonabelian Hodge theory, Homological mirror symmetry, perverse sheave, geometric Langlands program.
 These relations manifest themselves by the appearance of the following structures on $\JA$:

1) a sheaf of reductive Lie algebras,

2) (singular) Fano toric varieties whose hyperplane sections are (singular) Calabi-Yau varieties,

3) trivalent graphs.

This is an expository paper giving an account of most of the main properties of $\JA$
 uncovered in \cite{[R1]} and \cite{[R2]} . 
\end{abstract}
\tableofcontents
\section{Introduction}

 The Jacobian of a smooth complex projective curve is one of the most remarkable objects in algebraic geometry.
 It is well-known that the Jacobian $J(C)$ of a smooth complex projective curve $C$ is a principally polarized abelian variety whose dimension
equals the genus $g_C$ of the curve $C$. Its function theory reflects in a beautiful way a group structure of $J(C)$ as well as its relation to
 geometry of  $C$. This is  a very rich theory of theta-functions (see \cite{[M]}, for a beautiful introduction to the theory of 
 Jacobians of curves).
\\
\indent
For the purposes of our story we view $J(C)$ as the moduli space of line bundles with a fixed Chern class (=degree) on $C$. Take this degree
to be zero. Then upon a choice of a base point, call it $p_0$, the relation of $C$ and $J(C)$ is given by the Abel-Jacobi maps
\begin{equation}\label{AJ}
\alpha^d_{p_0} : C^{[d]} \longrightarrow J(C)
\end{equation}
where $C^{[d]} $ is the $d$-th symmetric power of $C$ and the map $\alpha^d_{p_0}$ is defined by sending $Z=p_1 + \cdots + p_d$ to $\OO_C (Z - dp_0)$. 
These maps, for various $d$, encapsulate all essential features of $C$.
For example, the fibre of $\alpha^d_{p_0}$ over a point of $J(C)$ corresponding to a line bundle $\cal L$ is the complete linear system
$\mid {\cal L} (d p_0) \mid$. This is the Abel-Jacobi theorem. When this is nonempty taking a divisor $D$ in $\mid {\cal L} (d p_0) \mid$
and  analyzing the differential of $\alpha^d_{p_0}$ at $D$ one recovers the Riemann-Roch theorem for $\OO_C (D)$.
\\
\indent
The theory of special divisors on $C$, the Brill-Noether theory, can be approached through a study of the images of the maps 
$\alpha^d_{p_0}$.
\\
\indent
The case $d=g_C -1$ , for $g_C \geq 2$, is of special importance since the image of  $\alpha^{g_C -1}_{p_0} $
 coincides, up to translation, with the theta-divisor 
$\Theta_C$
of $J(C)$, the zero-locus of the Riemann
theta-function. Then, from the classical theorem of Torelli, one knows that the pair $(J(C), \Theta_C )$ determines $C$ up to an isomorphism.
\\
\indent
This rapid account is enough to see that the Jacobian of a smooth projective curve is an important tool to study geometry of curves as well
as an incredibly rich and beautiful object to be studied on its own.
\\
\indent
When one turns to the higher dimensional projective varieties one quickly discovers that a comparable object does not exist.
In \cite{[R1]} and \cite{[R2]} we initiated a study of a new version of Jacobian for a smooth complex projective surface. 
The object that we call nonabelian Jacobian possesses several classical aspects:
\\
1) it relates to the rank 2 vector bundles on a surface in question as well as to its Hilbert schemes of zero dimensional subschemes - the
proper replacement of the symmetric power in (\ref{AJ});
\\
2) it carries a distinguished divisor that ``sees" a great deal of geometry of our surface thus allowing to address the problem of special
vector bundles on a surface as well as Torelli problem.

But a new phenomenon emerges: our Jacobian comes along with a sort of a period map. This feature leads to a natural appearance of
various new structures:
\\
3) a sheaf of reductive Lie algebras on the nonabelian Jacobian 
\\
4) (singular) Fano toric varieties whose hyperplane sections are (singular) Calabi-Yau varieties
\\
5) trivalent graphs. 

It should be clear that 3) could be envisaged as a nonabelian analogue of the (abelian) Lie algebra structure of the classical Jacobian.
This feature naturally relates geometry of surfaces with representation theory of reductive Lie algebras/groups. Thus one has powerful 
methods of the representation theory to address various problems in the theory of surfaces as well as the possibility to define new invariants of the representation theoretic origin
for surfaces and 
vector bundles on surfaces.

The Fano varieties in 4) emerge naturally as the parameter spaces of natural families
of Higgs structures related to the variational aspects of the period maps attached to our Jacobian. It could be viewed as a nonabelian Hodge theory
in a spirit of Simpson associated to the nonabelian Jacobian. This aspect of the theory should produce new invariants for surfaces and vector bundles
on them coming from toric Fano varieties and Calabi-Yau varieties.

The trivalent graphs in 5) come in as a convenient pictorial device to record the period map and its variational aspects. But we believe that it also
points to relations to knot theory,
moduli spaces of curves (via the mechanism of ribbon graphs) as well as Conformal field theory.

Being such a multifaceted object it seems that our nonabelian Jacobian is a ``correct" counterpart of the classical Jacobian and
deserves a serious study.

The main purpose of these notes is to give a concise and informal account of the results in \cite{[R1], [R2]}. Our exposition concentrates on the main
ideas and we largely bypass the technical issues and give no proofs. For those the interested reader is referred to the above cited works. 
\\
\\
\indent
The paper is divided in two parts. In {\bf Part I} we review all the foundational results of \cite{[R1]}. Starting from the definition in \S2
we continue with the discussion of all the features of the nonabelian Jacobian in the subsequent sections. In \S6 we discuss an example
of complete intersections, where everything can be computed explicitly. The part is ended, \S7, with a summary
 of the results discussed. 

The second part of the paper is entirely devoted to the representation theoretic aspects of our Jacobian and constitutes, with the exception of \S \ref{Langlands}, the contents of \cite{[R2]}. We tried to make it relatively self-contained.
So the reader interested only in this feature could go directly to {\bf Part II} and refer to the first part occasionally for notation and a  bare minimum of
explanations.  
\\
\\
\indent
This work has been started during the author's visit of the Hebrew University of Jerusalem, in October, 2010.
The author expresses his gratitude to this institution for financial support, hospitality and excellent working conditions. 
It is a pleasure to thank David Kazhdan for his interest in this work and for numerous stimulating discussions which helped to bring to
focus many aspects of this exposition.

\part{Nonabelian Jacobian of smooth projective surfaces}

 In our quest for a suitable version of Jacobian we are guided by two essential features of the classical Jacobian:
\\
\\
1) it is a parameter space of vector bundles (of rank 1) with fixed topology\label{guide}
\\
2) it relates to the symmetric product of a curve ( via the Abel-Jacobi maps in (\ref{AJ})).
\\
\\
Thus we are led to seek an object which parametrizes vector bundles on a surface together with certain data related to zero-dimensional
subschemes on a surface. The simplest candidate that comes to mind would be a vector bundle together with a global section vanishing on 
a zero-dimensional subscheme.  Since we are on a surface we are forced to consider vector bundles of rank 2. Thus we arrive to the conclusion that
a candidate for a nonabelian Jacobian of a surface $X$ should involve pairs $(\SE,e)$, where $\SE$ is a vector bundle of rank 2 on $X$ and
$e$ is its global section whose scheme of zeros is 0-dimensional. The construction and various properties of such a scheme of pairs is the contents of
\cite{[R1]}. This part discusses the essential aspects of that work.    
\section{Construction and properties}
{\bf 2.1. Scheme of pairs.} Let $X$ be a smooth complex projective surface. Fix a divisor $L$ on $X$ and an integer $d>0$.
To simplify the discussion we assume the divisor $L$ to be subject to the following properties:
$$
H^0 (\OO_X (-L)) =H^1 (\OO_X (-L))=0.
$$
This means that $L$ is assumed to be sufficiently `positive', e.g. ample or, more generally, numerically effective and big\footnote{i.e. $L^2 >0$.}.
Since we are aiming at geometric applications such assumptions are reasonable.

The object called a nonabelian Jacobian  $\JA$ of type $(L,d)$ is the universal scheme parametrizing pairs 
$(\SEE)$, where $\SE$ is a torsion free sheaf 
 of rank 2 with $det(\SE) =\OO_X (L)$, and
$d$, the degree of the second Chern class,  and $[e]$ is the homothety class of
a global section $e$  
of $\SE$, whose zero locus $Z_e$ is 0-dimensional.\footnote{the notion of zero-locus of a global section of a torsion free sheaf is ill defined, so to be precise,
instead of the homothety class of
a global section $e$  
of a torsion free sheaf $\SE$ we should speak about a monomorphism $\OO_X \stackrel{e}{\longrightarrow} \SE$; to say that $e$ has the zero-locus $Z_e$ of dimension $0$ means that the cokernel
 of $e$ is a torsion free sheaf; since that sheaf is of rank $1$ it has the form $\ID_Z (L)$, where $\ID_Z$ is the sheaf of ideals of a $0$-dimensional subscheme $Z$ of $X$;
 we {\it define} it to be $Z_e$, the zero-locus of (the homothety class of) $e$. We are thankful to P.Deligne for pointing out to us this subtlety.}   From this it follows that $\JA$ is related, on the one hand, to the stack 
${\bf M}_X (2,L,d)$ of torsion free sheaves on $X$ having rank 2, fixed determinant $\OO_X (L)$ and the second Chern class
$d$ and, on the other hand, to the Hilbert scheme $\XD$ of $0$-dimensional subschemes of $X$ having length $d$. 
This is recorded by the following diagram 
\begin{equation}\label{pi-h}
\xymatrix{
& {\JA} \ar[dl]_{\pi} \ar[dr]^{h}& \\
\XD & &{\bf M}_X (2,L,d) }
\end{equation}
The morphism $h$ takes the pair $(\SEE)$ to the object $[\SE]$ of ${\bf M}_X (2,L,d)$ corresponding to the sheaf $\SE$, while $\pi$ sends
$(\SEE)$ to the point $[Z_e]$ of the Hilbert scheme $\XD$ corresponding to the subscheme of zeros $Z_e$ of the section $e$.
\\
\indent
The fibre $h^{-1} ([\SE])$ of $h$ over $[\SE]$ is easily seen to be isomorphic to the Zariski open subset 
$U_{\SE}$ of $\PP(H^0 (\SE))$ parametrizing (homothety classes of) global sections of $\SE$ with 0-dimensional zero-locus.
This could be viewed as an analogue of the Abel-Jacobi theorem in our setting. One can rephrase this by saying that the morphism
$h$ defines a rational equivalence relation, denoted $\sim_{\scriptscriptstyle{AJ}}$ (`AJ' stands for Abel-Jacobi), such that the quotient
\begin{equation}\label{qAJ}
\JA /{\sim_{\scriptscriptstyle{AJ}}} = \mbox{is isomorphic to a substack of}\,\,{\bf M}_X (2,L,d). 
\end{equation}
\indent
Next we turn to the morphism $\pi$.
The fibre of $\pi$ over a point $[Z] \in \XD$ can be identified with torsion free sheaves $\SE$ fitting into the following short exact sequence
of sheaves on $X$
\begin{equation}\label{ext-seq}
\xymatrix{
0 \ar[r] & {\OO_X}  \ar[r] & {\SE} \ar[r] & {\ID_{Z} (L)} \ar[r]& 0 }
\end{equation}
where $\ID_Z$ is the ideal sheaf of the subscheme $Z$ in $X$. Such exact sequences are parametrized by the group of extensions
$Ext^1 (\ID_{Z} (L), \OO_X)$. In the sequel, to simplify the notation, this group will be denoted by
$Ext^1_Z$.
  Multiplying the extension class $\alpha \in Ext^1_Z$ corresponding to the sequence (\ref{ext-seq}) by a nonzero scalar does not change
the sheaf $\SE$ in the middle. This leads to the following identification
\begin{equation}\label{fib-pi}
\pi^{-1} ([Z]) \cong \PP(Ext^1_Z)
\end{equation}
Thus the morphism $\pi$ defines another rational equivalence relation on $\JA$, denoted by $\sim_{\scriptscriptstyle{S}}$ (`S' is for Serre 
since the extension sequence in (\ref{ext-seq}) is an instance of the Serre construction associating a torsion free sheaf with a codimension
2 subscheme on a smooth projective variety, see \cite{[O-S-S]}).  Taking the quotient of $\JA$ by this equivalence relation one obtains a subscheme of
the Hilbert scheme. The points $[Z]$ belonging to this subscheme acquire now ``internal" moduli - the projective space
$\PP(\EZ)$ in (\ref{fib-pi}), which is nontrivial as soon as $dim(\EZ) \geq 2$. It is by exploiting these nontrivial ``internal" moduli that new structures of
$\JA$ are uncovered. It also should be remarked that for a 0-dimensional subscheme to have ``internal" moduli is an intrinsic feature of higher dimensional
geometry because, 
 if one applies the guiding principles 1)-2) (see the top of page \pageref{guide}) to curves, then one deals with pairs
$({\cal L}, [e])$ consisting of a line bundle $\cal L$ on a curve and the homothety class $[e]$ of a nonzero global section $e$ of $\cal L$.
Then a given 0-dimensional subscheme $Z$ on a curve $C$ determines a unique line bundle, namely $\OO_C (Z)$, having a section vanishing on $Z$.

In what follows we explore $\JA$ as a scheme\footnote{formally $\JA$ is defined as ${\bf Proj} ({\cal S})$ of a certain coherent sheaf ${\cal S}$ on $\XD$ (see \RI, \S1.1). 
So it comes equipped 
with the structure morphism $\pi: \JA \longrightarrow \XD$ as well as a canonical choice of an invertible sheaf $\OO_{\JA} (1)$ on $\JA$ such that
the direct image $\pi_{\ast} (\OO_{\JA} (1)) ={\cal S}$.}\label{proj} over $\XD$ via the morphism $\pi$. From this point of view closed points of $\JA$ are viewed as pairs
$\ZA$, where $[Z] \in Im(\pi)$ and $[\alpha]$ is a point in $\pi^{-1} ([Z]) = \PP(\EZ)$. To go back and forth between two notation,
$(\SEE)$ and $\ZA$, is facilitated by the extension sequence in (\ref{ext-seq}). In the sequel we switch between these two ways to denote
points of $\JA$ without further comments. 
\\
\\
{\bf 2.2. Stratification of $\XD$ and $\JA$.}
The morphism $\pi :\JA \longrightarrow \XD$ induces the stratification of the Hilbert scheme
$\XD$ according to the dimension of fibres of this morphism.  
Namely, for every integer $r\geq 0$ we set
$$
\Gamma^{r}_{d} (L )= \{ [Z] \in \XD \mid dim(\pi^{-1} ([Z])) = dim(\PP (\EZ)) \geq r \}
$$
This gives the following stratification of $\XD$
\begin{equation}\label{strat-Gam}
 X^{[d]} \supset \,\Gamma^{0}_{d} (L )\,\supset\,\Gamma^{1}_{d} (L )\,\supset 
  \dots \supset\,\Gamma^{r}_{d} (L )\,\supset \dots
 \end{equation}
The sets $\Gamma^{r}_{d} (L )$ are in fact the degeneracy loci of a morphism between vector bundles on $\XD$ and
hence they carry a natural structure of closed subvarieties of $\XD$ (see \cite{[R2]}, \S1.2). 

 The natural projection $\pi: \JA \longrightarrow \XD$
 induces the stratification 
 \begin{equation}\label{strat-Jac}
 \JA = \JAA^0 \supset \JAA^1 \supset \dots \supset \JAA^r \supset \dots\
 \end{equation}
 where $\JAA^r = \pi^{-1}({\Gamma^{r}_{d} (L )})$. In particular,
each stratum $\JAA^r$ is a closed subscheme of $\JA$ and the locally closed stratum 
 $\stackrel{\circ}{\JAA^r} = \JAA^r \setminus \JAA^{r+1}$, when non-empty, is a 
 $\PP^r$-bundle over 
 $\stackrel{\circ}{\Gamma^{r}_{d}} (L) =\Gamma^{r}_{d} (L ) \setminus
 \Gamma^{r+1}_{d} (L ) $.
 Of course, we get something new only for $r\geq1$. This will be assumed
 for the rest of the paper.
\\
\\
{\bf 2.3. The nonabelian Theta-divisor ${\bf\Theta} (X;L,d)$.} 
Similar to the classical Jacobian the scheme $\JA$ comes equipped      
 with a distinguished divisor. Some care has to be taken to define it (see \cite{[R1]}, 1.2.) but as a set of points
${\bf\Theta} (X;L,d)$ is easy to grasp - it parametrizes pairs
$(\SEE)$, where $\SE$ is {\it not} locally free. In particular, the fibre
${\bf \Theta}_Z$ of ${\bf\Theta} (X;L,d)$ over $[Z] \in \XD$ is a hypersurface of degree $d$ in 
$\PP(\EZ)$. More precisely, one can show that, set-theoretically,
${\bf\Theta}_Z$ is the union of hyperplanes 
$H_z$ in $\PP(\EZ)$, where $z$ runs through the set of closed points of $Z$.
Thus the divisor ${\bf\Theta} (X;L,d)$ captures the geometry of zero-dimensional subschemes of $X$ parametrized by the
underlying points of the Hilbert scheme $\XD$.

 At this stage the pair $(\JA,{\bf \Theta}(X;L,d))$ 
 is a rather precise analogue of the classical
 Jacobian. However, there is a new feature: $\JA$ carries a
 variation of Hodge-like structure.

\section{ Variation of Hodge-like structure on $\JA$}
By a Variation of Hodge-like structure we mean a decreasing filtration of 
 sheaves equipped with a ``derivative"
  which shifts the index of the filtration
 at most by $1$ (an analogue of Griffiths transversality property for
 the Infinitesimal variation of Hodge structure, \cite{[G]}).
It turns out that certain basic sheaves on $\JA$ come along with such structure.
This seems to us a qualitatively new feature of our Jacobian (with respect
to the classical one). It reflects at the same time the fact that we are 
in dimension $>1$ and are dealing with higher rank bundles.
\\
\\
{\bf 3.1. The sheaf $\FT$ and its filtration $\HT_{\bullet}$.} To begin our construction we recall that
over the Hilbert scheme $\XD$ of closed zero-dimensional subschemes of $X$ having length $d$
there is the universal scheme $\ZD$ of such subschemes
\begin{equation}\label{uc}
\xymatrix{
&\ZD \ar[dl]_{p_1} \ar[dr]^{p_2} \ar@{^{(}{-}{>}}[r]& X \times \XD\\
X& &\XD &}
\end{equation}
where $p_i (i=1,2)$ is the restriction to $\ZD$ of the projections $pr_i,\,\,(i=1,2)$ of the Cartesian product $X \times \XD$
onto the corresponding factor. For a point $\xi \in \XD$ the fibre $p^{\ast}_2 (\xi)$ is isomorphic via $p_1$ with the subscheme $Z_{\xi}$ of $X$ corresponding to $\xi$, i.e.
\begin{equation}\label{fc}
 Z_{\xi} = p_{1 \ast} (p^{\ast}_2 (\xi) )\,.
\end{equation}
In the sequel we often make no distinction between $Z_{\xi}$  
and the fibre  
$
p^{\ast}_2 (\xi) 
$
itself. If $Z$ is a closed subscheme of dimension zero and length $d$, then $[Z]$ will denote the corresponding point in the Hilbert scheme
$
\XD
$.

One of the basic sheaves on $\XD$ is the direct image
$p_{2 \ast} \OO_{\ZD}$ of the structure sheaf $\OO_{\ZD}$ of $\ZD$. This is a locally free sheaf of rank $d$, whose fibre over a point 
$[Z] \in \XD$ is the space
$H^0 (\OO_Z)$ of complex valued functions on $Z$. We will be interested  in its pullback 
\begin{equation}\label{tildF}
\FT = \PIPO
\end{equation}
to $\JA$. This is a natural locally free sheaf on $\JA$. Its fibre $\FT_{\ZA}$ at a point\footnote{Here we switch to the alternative, equivalent, notation
of points in $\JA$ discussed in the end of \S2.1.} $\ZA$ of $\JA$ 
is again $ H^0 (\OO_{Z})$. However, the points of $\ZA$ of $\JA$ lying over $[Z]$ give an
additional structure to the space $H^0 (\OO_{Z})$. This is recorded in the following 
\begin{lem}\label{Htild}
 Let $[Z]$ be a point in $\Gamma^{r}_{d} (L)$, for $r\geq0$, and let $\ZA$ be a point of $\JAA^r$ lying over $[Z]$. 
Then there is a distinguished linear subspace
$\HT \ZA$ of $H^0 (\OO_{Z})$ having the following properties:
\begin{enumerate}
\item[1)]
$dim (\HT \ZA) \geq r+1$,
\item[2)]
$\HT \ZA$ contains the subspace of constant functions on $Z$,
\item[3)]
If $\ZA =(\SEE)$ is such that $\SE$ is locally free, then $\HT \ZA$ is naturally isomorphic\footnote{an isomorphism is unique up to a non-zero scalar.} to the group of
extensions $\EZ$. Moreover, under this isomorphism the subspace of constant functions on $Z$ in $\HT \ZA$ goes over to the line $\CC\{\alpha\}$ in $\EZ$.
\end{enumerate}
\end{lem}
 
Varying $\ZA$  in $\JAA^r$, the spaces 
$\HT \ZA$ fit together to form the subsheaf 
$\HT$ of $\FT$. The sheaf $\HT$ is the first step of a distinguished filtration of $\FT$ whose construction and basic properties
are summarized in the proposition below. 
\begin{pro-defi}\label{filtrn}
On every stratum $\JAA^r, (r\geq1)$
 the sheaf of rings $\FT=\PIPO$ admits a distinguished filtration
 \begin{equation}\label{H-filt}
 0={\HT_{0}} \subset {\HT_{-1}} \subset \dots
 \subset \FT
   \end{equation}
subject to the following properties.
\\
a) ${\HT_{-1}} = \HT$ and it contains the structure sheaf $\OO_{\JAA^r}$.
\\
b) The multiplication in $\FT$ induces
the morphisms
$$
{\bf m_k}:S^k \HT \longrightarrow \FT\,,
$$
for every $k\geq 1$. The sheaf $\HT_{-k}$ is defined to be the image of ${\bf m_k}$.
In particular, one obtains the multiplication morphisms
\begin{equation}\label{formalGr-H}
\HT \otimes \HT_{-k} \longrightarrow \HT_{-k-1}\,,
\end{equation}
for every $k\geq 1$. 
\end{pro-defi}

The filtration $\HT_{\bullet}$ in (\ref{H-filt}) allows to partition irreducible components of $\JAA^r$ according
to the ranks of the constituents of that filtration. More precisely, fix a stratum
$\JAA^r$ so that that for a general point $(\SEE) \in \JAA^r$ the sheaf $\SE$ is locally free and the locally closed stratum
$\stackrel{\circ}{\JAA^r} = \JAA^r \setminus \JAA^{r+1}$ is non-empty. This is a $\PP^r$-bundle over
the stratum
$\stackrel{\circ}{\Gamma^{r}_{d}} (L) =\Gamma^{r}_{d} (L ) \setminus
 \Gamma^{r+1}_{d} (L ) $
of the stratification of $\XD$ in (\ref{strat-Gam}).
In particular, if we let $C^r (L,d)$ to be the set of irreducible components of
 $\Gamma^{r}_{d} (L)$, then 
$$
\{\JAA_{\Gamma} =\pi^{-1} (\Gamma) \mid \Gamma \in C^r (L,d) \}
$$
is the set of the irreducible components of $\JAA^r$.
 
On every irreducible component $\JAA_{\Gamma}$ of $\JAA^r$ the sheaves $\HT_{-i} \otimes \OO_{\JAA_{\Gamma}}$ are non-zero and torsion free,
 for every $i\geq 1$. 
 So their ranks are well-defined. Set
\begin{equation}\label{h-Gam}
h^{i-1}_{\Gamma} = rk(\HT_{-i} \otimes \OO_{\JAA_{\Gamma}}) - rk(\HT_{-i+1} \otimes \OO_{\JAA_{\Gamma}})\,,
\end{equation}
for every $i\geq 1$.

 Denote by $l_{\GA}$ the largest index $i$ for which
$h^{i-1}_{\GA} \neq 0$
and call it {\it the length of the filtration} of
$\HT_{\bullet} \otimes \OO_{\JG}$.
Thus on 
$\JG$
the filtration (\ref{H-filt}) stabilizes at 
$\HT_{-l_{\GA}} \otimes \OO_{\JG}$.
We also agree to assign to
$\FT\otimes \OO_{\JG}$
the index 
$-(l_{\GA} +1)$
and use the notation
\begin{equation}\label{lterm}
\FT\otimes \OO_{\JG} = \HT_{-l_{\GA} -1}\,.
\end{equation}
So the filtration (\ref{H-filt}) restricted to $\JG$ has the following form
\begin{equation}\label{filtHT-JG}
0=\HT_0 \otimes \OO_{\JG} \subset \HT_{-1} \otimes \OO_{\JG} \subset \ldots \subset \HT_{-l_{\GA} +1} \otimes \OO_{\JG}
 \subset  \HT_{-l_{\GA}} \otimes \OO_{\JG} \subset   \HT_{-l_{\GA}-1} = \FT \otimes \OO_{\JG}
\end{equation}
and it attaches to a component $\Gamma \in \CS$ the vector
\begin{equation}\label{Hilb-Gam}
h_{\GA} =(h^0_{\GA},\ldots, h^{l_{\GA} -1} _{\GA},h^{l_{\GA} } _{\GA})
\end{equation}
which will be called {\it the Hilbert vector of $\GA$}.
\begin{lem}\label{hg=c}
The Hilbert vector $h_{\GA}$ and its components
$h^i_{\GA},\,\,(i=0,\ldots, l_{\GA}),$
have the following properties:
\begin{enumerate}
\item[1)]
$h_{\GA}$ is a composition of $d$, i.e.
$$
\sum^{l_{\GA}}_{i=0} h^i_{\GA} = d\,.
$$
\item[2)]
$h^0_{\GA}= rk(\HT_{-1} \otimes \OO_{\JG}) =rk(\HT \otimes \OO_{\JG}) = r+1$.
\item[3)]
$h^i_{\GA} >0$, for $i=0,\ldots, l_{\GA} -1$, and
$h^{l_{\GA}}_{\GA} \geq 0$.
\end{enumerate}
\end{lem}

Let $C(d)$ be the set of compositions of $d$. From Lemma \ref{hg=c} it follows
that the assignment of the Hilbert vector
$h_{\GA}$
to the components $\GA$ in $\CS$ gives a map
\begin{equation}\label{CS-Hv}
\xymatrix@1{
{h(L,d,r): \CS} \ar[r]& C(d)}\,.
\end{equation}

The algebro-geometric meaning of the filtration (\ref{H-filt}) is easy to understand:
let ${\HT}([Z],[\alpha])$ be the fibre of $\HT$ at a closed point
$([Z],[\alpha])$ of $\JAA^r$. This is a subspace of $\HO Z )$. So it can be viewed as a linear
system on $Z$. Moreover, the linear system is base point free since it contains the 
subspace $\CC\{1_Z\}$ of the constant functions\footnote{$1_Z$ here and further on denotes the constant function of value $1$ on $Z$.} on $Z$. So 
${\HT}([Z],[\alpha])$ defines a morphism
\begin{equation}\label{kappaZ}
\kappa_{([Z],[\alpha])} : Z \longrightarrow \PP({\HT([Z],[\alpha])}^{\ast})\,.
\end{equation}
Thus the filtration (\ref{H-filt}) captures the geometry of the morphism
$\kappa_{([Z],[\alpha])}$ and the Hilbert vector $h_{\Gamma}$ encodes
 the Hilbert function of the image of 
$\kappa_{([Z],[\alpha])}$, for all $\ZA$ varying in the complement of the singularity loci\footnote{the singularity locus of a sheaf ${\cal G}$ on a scheme $Y$ is the subscheme of $Y$, where  ${\cal G}$
is not locally free.} of the sheaves
$\HT_{-i} \otimes \OO_{\JG}$, for $i=1,\ldots,l_{\GA}$. This complement is a non-empty Zariski open subset of 
$\JG$ which will be denoted $\JAA^{\prime}_{\GA}$. 

Set
\begin{equation}\label{G(0)}
\GA^{(0)} = \pi (\JG^{\prime})
\end{equation}
to be the image of $\JG^{\prime}$ under the projection $\pi$ in (\ref{pi-h}). This is a Zariski open subset of
$\GA$. The following result relates $\JAA^{\prime}_{\GA}$ to the complement of the theta-divisor $\TE(X;L,d)$ in $\pi^{-1} ( \GA^{(0)})$. 
\begin{lem}\label{lem-rk}
Let 
$ \TE_{\GA^{(0)}} =\TE(X;L,d) \cap \pi^{-1} ( \GA^{(0)})$
be the theta-divisor over 
$\GA^{(0)}$
and let
$$
\JG^{(0)} = \pi^{-1} ( \GA^{(0)}) \setminus  \TE_{\GA^{(0)}}
$$
be its complement in
$\pi^{-1} ( \GA^{(0)})$.
Then
$
\JG^{(0)} \subset \JG^{\prime}
$.
\end{lem}
\vspace{0.5cm}

{\bf 3.2. Orthogonal decomposition of $\FT$.}\label{sec-ord}
The filtration $\HT_{\bullet}$ in (\ref{H-filt}) acquires more structure over the points of
$\GG(L) \subset \XD$ corresponding to the reduced subschemes of $X$.
Let 
$Conf_d (X)$
be the locus of the Hilbert scheme 
$\XD$
parametrizing the subschemes of $d$ distinct points of $X$. This is a Zariski open subset of
$\XD$ since it can be described as the complement of the branching divisor of the ramified covering
$$
\xymatrix@1{
{p_2 :\ZD} \ar[r]& {\XD} }
$$
in (\ref{uc}). The subschemes $Z$ of $X$ with
$[Z] \in \CO$ will be called configurations (of $d$ points) on $X$.
\\
\indent
We are interested in the irreducible components 
$\GA \in \CS$
having a non-empty intersection with
$\CO$.
\begin{defi}\label{csa}
A component 
$\GA \in \CS$
is called admissible if
$$
\GA \cap \CO \neq \emptyset\,.
$$
The set of admissible components in $\CS$ will be denoted by
$C^r_{adm} (L,d)$.
\end{defi}
For a subset $Y$ in
$\XD$
we denote by
$Y_{conf}$
the intersection
$Y \cap \CO$.
In particular, for  
$\GA \in \CS$
the subset $\GA_{conf}$ is Zariski open in $\GA$ and it is non-empty if and only if 
$\GA \in \CSA$ .
The subset $\GA_{conf}$ will be called the configuration subset of $\GA$.
We will now explain why configurations are important in our constructions (see \RI, \S2, for details).
\\
\\
\indent
The sheaf $p_{2\ast} \OO_{\ZD}$ as well as its pullback $\FT =\PIPO$ admit the trace morphism
\begin{equation}\label{tr}
\xymatrix@1{
{Tr: \FF=p_{2\ast} \OO_{\ZD}} \ar[r]& {\OO_{\XD}},& & {\tilde{Tr}: \FT} \ar[r]& {\OO_{\JA}}. }
\end{equation}
It can be used to define the bilinear, symmetric pairing
$\QB$ (resp. ${\bf \tilde{q}}$) on $\FF$ (resp. $\FT$)
\begin{equation}\label{q}
\QB(f,g) =Tr(fg)\,\,({\rm{resp.}}\,\, {\bf \tilde{q}}(f,g) = \tilde{Tr} (fg))\,,
\end{equation}
for every pair $(f,g)$ of local sections of $\FF$ (resp. $\FT$).

This pairing is non-degenerate precisely over
$\CO$. Using it we obtain a natural splitting of the filtration
$\HT_{\bullet}$ in (\ref{H-filt})
on a certain Zariski open subset of
$\JAA_{\GA^{(0)}_{conf}} =\pi^{-1} (\GA^{(0)}_{conf})$,
for every admissible component
$\GA \in \CSA$ (recall $\GA^{(0)}$ is a Zariski open subset of $\GA$ defined in (\ref{G(0)})).
 More precisely, set
\begin{equation}\label{ort-c}
\FI^i =\big{(} \HT_{-i} \otimes \OO_{\JAA_{\GA^{(0)}_{conf}}} \big{)}^{\perp}
\end{equation}
to be the subsheaf of
$\FT  \otimes \OO_{\JAA_{\GA^{(0)}_{conf}}}$
orthogonal to 
$ \HT_{-i} \otimes \OO_{\JAA_{\GA^{(0)}_{conf}}} $
with respect to the quadratic form
$\QT$ in (\ref{q}). It was shown in \RI, Corollary 2.4, that there exists a non-empty Zariski open subset
$\JAB_{\GA}$ of $\JG$ subject to the following properties:
\begin{enumerate}
\item[(a)]
the open set 
$\JAB_{\GA}$ lies over 
$ \GA^{(0)}_{conf}$, i.e. the morphism $\pi$ in (\ref{pi-h}) restricted to 
$\JAB_{\GA}$ gives the surjective morphism
\begin{equation}\label{pi-JAB}
\xymatrix@1{
{\pi: \JAB_{\GA}} \ar[r]& {\GA^{(0)}_{conf}} }\,.
\end{equation}
\item[(b)]
$\JAB_{\GA}$ lies in the complement of the theta-divisor, i.e.
\begin{equation}\label{JAB-inc}
\JAB_{\GA} \subset \JG^{(0)},
\end{equation}
 where $ \JG^{(0)} $ is as in Lemma \ref{lem-rk}.  
\item[(c)]
$\FT $ restricted to $\JAB_{\GA}$ admits the orthogonal direct sum decomposition
\begin{equation}\label{FT-ors}
\FT \otimes \OO_{\JAB_{\GA}} =\HT_{-i} \otimes \OO_{\JAB_{\GA}} \oplus \FI^i  \otimes \OO_{\JAB_{\GA}},
\end{equation}
for every $i=0, 1, \ldots, l_{\GA} +1$.
\end{enumerate}
This gives rise to the filtration
\begin{equation}\label{filtF}
\FT \otimes \OO_{\JAB_{\GA}} =\FI^0  \otimes \OO_{\JAB_{\GA}} \supset \FI^1  \otimes \OO_{\JAB_{\GA}} \supset \ldots \supset 
\FI^{ l_{\GA}} \otimes \OO_{\JAB_{\GA}} \supset \FI^{ l_{\GA} +1} =0.
\end{equation}
Putting together the filtrations $\HT_{\bullet} \otimes \OO_{\JAB_{\GA}}$ and $\FI^{\bullet}  \otimes \OO_{\JAB_{\GA}}$
we define the subsheaves
\begin{equation}\label{H}
\HH^{i-1} =\big{(} \HT_{-i} \otimes \OO_{\JAB_{\GA}} \big{)}\cap \big{(} \FI^{i-1}  \otimes \OO_{\JAB_{\GA}}\big{)},
\,\,{\rm{for}} \,\,i=1,\ldots,  l_{\GA} +1.
\end{equation}
This definition together with (\ref{FT-ors}) yield the following decomposition of 
$\FT  \otimes \OO_{\JAB_{\GA}}$
into the orthogonal sum
\begin{equation}\label{ordFT}
\FT  \otimes \OO_{\JAB_{\GA}} = \bigoplus^{ l_{\GA}}_{p=0} \HH^p .
\end{equation}
\begin{rem}\label{rkH}
\begin{enumerate}
\item[1)]
Observe that the ranks of the summands
$\HH^p$'s in (\ref{ordFT}) form the Hilbert vector
$h_{\GA}$ as defined in Lemma \ref{hg=c}, i.e.
\begin{equation}\label{rkHp}
rk(\HH^p)= h^p_{\GA}.
\end{equation}
This follows from the definition of 
$ h^p_{\GA}$ in (\ref{h-Gam}), the inclusion (\ref{JAB-inc})
and the orthogonal decomposition of $\HT_{-i} \otimes \OO_{\JAB_{\GA}}$
\begin{equation}\label{ordH-i}
\HT_{-i} \otimes \OO_{\JAB_{\GA}} = \bigoplus^{ i-1}_{p=0} \HH^p .
\end{equation}
In particular, for $i=1$ one obtains
\BEN\label{ordH-1}
\HT =\HT_{-1} =\HH^0.
\EEN
\item[2)]
From the orthogonal decomposition (\ref{ordFT}) it follows that the subsheaves 
$\FI^i$ of the filtration $\FI^{\bullet}$ in (\ref{filtF}) admit the following orthogonal decomposition
\BEN\label{ordF-i}
\FI^i = \bigoplus^{\LG}_{p=i} \HH^p .
\EEN
\end{enumerate}
\end{rem}

For $\ZA \in \JAB_{\GA}$ the decomposition (\ref{ordFT}) says that the only non-trivial de Rham cohomology group
$H^0 (\OO_Z)$ of $Z$ admits the direct sum decomposition
\begin{equation}\label{ord-ZA}
H^0 (\OO_Z) =\bigoplus^{ l_{\GA}}_{p=0} \HH^p \ZA ,
\end{equation}
where $\HH^p \ZA$ is the fibre of $\HH^p$ at $\ZA$. This could be viewed as a sort of Hodge decomposition for
(smooth) $0$-dimensional subschemes of $X$. More precisely, the above decomposition of $H^0 (\OO_Z)$ becomes visible in the presence of an extension class
 $[\alpha]$ in $\JAB_Z$, the fibre of $\JAB_{\GA}$ over $[Z]$. So such extensions could be envisaged as non-classical K\"ahler structures of $Z$ inducing
the Hodge-like decomposition of $H^0 (\OO_Z)$ in (\ref{ord-ZA}). Furthermore, the direct sum
 in (\ref{ordFT}) implies that this decomposition varies holomorphically as $\ZA$ moves in $\JAB_{\GA}$. Thus the situation resembles the
usual variation of Hodge structure and one can associate with $\JAB_{\GA}$ the corresponding period map.\footnote{In fact there are two natural choices in the situation at hand.
Both are treated in \cite{[R2]},\S4. In \S3.3, we discuss only one of them. The reader should be aware that that the period map denoted in \S3.3 by $p_{\GA}$ 
corresponds to ${}^{op} p_{\GA}$ in \cite{[R2]},\S4.}

Before we go on to discuss our period map, let us close this subsection by giving the dual version of the filtration\label{alt}
$\HT_{\bullet}$ in (\ref{filtHT-JG}) which coincides with the filtration
$\FI^{\bullet}$ in (\ref{filtF}) once restricted to $\JAB_{\GA}$. However, it has a virtue of being more geometric.\footnote{the subsequent discussion  is not used in the paper with
the exception of \S9.2, {\bf Example 1}, so the reader can skip it and continue to \S3.3.}

The starting point of the dual construction is another natural sheaf on $\XD$
which takes account of the divisor $L$. Namely, we consider the sheaf
\begin{equation}\label{def-FL}
\FF(L) = p_{2\ast} \big{(} p^{\ast}_1 \OO_X (L+K_X) \big{)}\,,
\end{equation}
where $p_i (i=1,2)$ are as in (\ref{uc}).  Taking its pullback via $\pi$ in 
(\ref{pi-JAB}) we obtain the sheaf
\begin{equation}\label{FTL}
\FT(L) = \pi^{\ast} \big{(} \FF(L) \big{)}\,.
\end{equation}
In \RI, \S1.3, it was shown that there is a natural morphism
\begin{equation}\label{mor-R}
\xymatrix@1{
{{\bf R^r}: \FT \otimes \OO_{\JAA^r}} \ar[r]& H^0(L+K_X)^{\ast} \otimes \OO_{\JAA^r} (1) }\,,
\end{equation}
where 
$\OO_{\JAA^r} (1) $
is the restriction to $\JAA^r$ of the tautological invertible sheaf
$ \OO_{\JA} (1)$  (see the footnote on page \pageref{proj} for notation).
In particular, the subsheaf 
$\HT$ which we encountered in Proposition-Definition \ref{filtrn}, a), is defined in \RI, (1.21), as the kernel of 
${\bf R^r}$ and we have morphisms
\begin{equation}\label{mor-RT-i}
\xymatrix@1{
{{\bf  \tilde{R}^r_i}: S^i \HT} \ar[r]^{m_i}&{\FT \otimes \OO_{\JAA^r}} \ar[r]^(0.35){\bf R^r}& H^0(L+K_X)^{\ast} \otimes \OO_{\JAA^r} (1) }\,,
\end{equation}
where $m_i$ is as in Proposition-Definition \ref{filtrn}, b). 

Dualizing and tensoring with 
$ \OO_{\JAA^r} (1)$
yields
$$
\xymatrix@1{
H^0(L+K_X) \otimes \OO_{\JAA^r}  \ar[r] & {\big{(} S^i \HT \big{)}^{\ast} \otimes \OO_{\JAA^r} (1) }\,. }
$$
Setting $\FIT_{i}$ to be the kernel of this morphism we obtain the following filtration
$$
H^0(L+K_X) \otimes \OO_{\JAA^r}= \FIT_1 \supset \FIT_2 \supset  \ldots \supset \FIT_i \supset \FIT_{i+1} \supset \ldots
$$
Each 
$\FIT_i $
contains the sheaf 
$
\pi^{\ast} \big{(} pr_{2\ast} \big{(} \JJ_{\ZD} \otimes pr^{\ast}_1 \OO_X (L+K_X) \big{)} \big{)}
$,
where 
$ \JJ_{\ZD}$
is the sheaf of ideals of the universal subscheme $\ZD$ in $X \times \XD$ (see (\ref{uc}) for notation) 
and $pr_j (j=1,2)$ are the projections of
$X \times \XD$
onto the corresponding factor
(this inclusion is proved in \RI, Proposition 1.6). Factoring out by 
$
\pi^{\ast} \big{(} pr_{2\ast} \big{(} \JJ_{\ZD} \otimes pr^{\ast}_1 \OO_X (L+K_X) \big{)} \big{)}
$
 one obtains the following filtration of 
$\FT(L)$:
\begin{equation}\label{filtFd}
\FT(L) \otimes \OO_{\JAA^r} =\FI_0   \supset \FI_1   \supset \ldots \supset \FI_i \supset \FI_{ i +1} \supset \ldots
\end{equation}
where
\begin{equation}\label{F-i}
 \FI_i  =\FIT_i  \big{/}  \pi^{\ast} \big{(} pr_{2\ast} \big{(} \JJ_{\ZD} \otimes pr^{\ast}_1 \OO_X (L+K_X) \big{)} \big{)}\,.
\end{equation}
To relate this filtration to the one in (\ref{filtF}) one observes that there is a natural morphism
\begin{equation}\label{mor-RT}
\xymatrix@1{
{{\bf \tilde{R}}: \FT  \otimes \OO_{\JAA^r}} \ar[r]& {(\FT (L))^{\ast} \otimes \OO_{\JAA^r} (1)} }
\end{equation}
(see \RI, (1.27) and (1.19) for details). Furthermore, this morphism is an isomorphism precisely on the complement of the theta-divisor
$\TE(X;L,d)$,  since the latter is defined by the vanishing of the determinant of 
${\bf \tilde{R}}$ 
(see the formula for $\TE(X;L,d)$ below \RI, (1.19)).
Taking the dual of
${\bf \tilde{R}}$ 
we obtain a natural identification of
$\FT (L) \otimes \OO_{\JAA^r} (-1) $
with
$ \FT^{\ast}  \otimes \OO_{\JAA^r}$
on the complement of the theta divisor in $\JAA^r$.
Restricting further to the configurations and using the self-duality of
$\FT$
over $\CO$ (provided by the quadratic form $\QT$ in (\ref{q}))
we obtain a natural identification of
$\FT (L) \otimes \OO_{\JAA^r} (-1) $
and
$ \FT  \otimes \OO_{\JAA^r}$
over the complement
$$
\pi^{-1} (\GAC) \setminus \TE(X;L,d)\,,
$$
for every component
$\GA \in \CSA$. In particular, this identification holds on 
$\JAB_{\GA}$ in view of the inclusion in (\ref{JAB-inc}).

With the above identification in hand, we can transfer the filtration
$\FI_{\bullet} \otimes \OO_{\JAB_{\GA}} (-1)$ of 
$\FT(L) \otimes \OO_{\JAB_{\GA}} (-1)$
in (\ref{filtFd}) (twisted by  $\OO_{\JAB_{\GA}} (-1)$)
to a filtration of
$\FT\otimes \OO_{\JAB_{\GA}}$.
 The point is that the resulting filtration is the filtration
$\FI^{\bullet}$ defined previously in (\ref{filtF}) via orthogonality
(see more detailed discussion in \RI, \S2).

By definition the filtration 
$\HT_{\bullet}$
is related to geometry of the morphisms
$\kappa \ZA$ in (\ref{kappaZ})
associated to the linear systems
$\left| \HT\ZA \right|$
on $Z$, as $[Z]$ varies through the points of the admissible components
$\GA$.
On the other hand the filtration $\FI^{\bullet}$ in (\ref{filtF}), in view
of its identification with
$\FI_{\bullet} \otimes \OO_{\JAB_{\GA}} (-1)$,
reflects geometric properties of the subschemes $Z$
(parametrized by $\GA$) with respect to the adjoint linear system
$\left| L+K_X \right|$ on $X$.
Thus one can say that the orthogonal decomposition
(\ref{ordFT}) contains information about the geometry of the subschemes
$Z$ with respect to {\it both} linear systems. 
\\
\\
{\bf 3.3.  Period map associated to $\JAB_{\GA}$.} We begin by defining the target for our period map.
This will be a certain scheme of partial flags determined by the Hilbert vector
$h_{\GA}=(h^0_{\GA},\ldots, h^{l_{\GA} -1} _{\GA},h^{l_{\GA} } _{\GA})$ in (\ref{Hilb-Gam}).
More precisely, set 
\begin{equation}\label{GAB}
\GA^{(0)}_{conf} = \GAB\,\,\,\,and \,\,\,\, \FF_{\GAB} =p_{2\ast} \OO_{\ZD} \otimes \OO_{\GAB}
\end{equation}
and
consider the scheme
$\mbox{\BM${\FL}$}_{h_{\GA}}$ of relative partial flags of type
$h_{\GA}$ in $ \FF_{\GAB}$, i.e.
$\mbox{\BM${\FL}$}_{h_{\GA}}$ is the scheme over $\GAB$ with the structure morphism
\BEN\label{fl}
Fl_{\GA} : \FLA \longrightarrow \GAB
\EEN
such that the fibre 
$\FLA ([Z])$ over a closed point $[Z] \in \GAB$ is the variety of partial flags of type
$h_{\GA}$ in the vector space
$\FF_{\GAB} ([Z]) = \HO Z)$, the fibre of $\FF_{\GAB}$ at $[Z]$. Hence the set of closed points of $\FLA ([Z])$ can be described as follows
\BEN\label{Fla-Z}
\FLA([Z]) =\left\{\left.[F] =[\scriptstyle{\HO{Z}) =F^0 \supset F^1 \supset \ldots \supset F^{\LG} \supset F^{\LG+1} =0}] \right| 
\scriptstyle{dim(F^{\LG  -p} /F^{\LG -p+1}) =h^{p}_{\GA},\,\,{\rm{for}}\,\, 0\leq p \leq \LG }\right\}\,.
\EEN
To relate $\JAB_{\GA}$ with $\mbox{\BM${\FL}$}_{h_{\GA}}$, we use the fact that the pullback
$\FT \otimes \OO_{\JAB_{\GA}}=\pi^{\ast} \FF_{\GAB}$ comes with the distinguished filtration $\HT_{\bullet} \otimes \OO_{\JAB_{\GA}}$, the restriction of the filtration (\ref{filtHT-JG}) to $\JAB_{\GA}$. 
Furthermore, by definition of 
$\JAB_{\GA}$ all sheaves of $\HT_{\bullet} \otimes \OO_{\JAB_{\GA}}$ are locally free and their successive quotients
\BEN\label{quot-HH}
 \HT_{-p} \otimes \OO_{\JAB_{\GA}}  / \HT_{-p+1} \otimes \OO_{\JAB_{\GA}} \cong \HH^{p-1}
\EEN
have ranks $h^p_{\GA}$, for $p=1,\ldots,\LG+1$ (the isomorphism in (\ref{quot-HH}) follows from the orthogonal decompositions in (\ref{ordH-i})). 
By the universality of 
$\FLA$
we obtain  the unique morphism
\BEN\label{p}
p_{\GA} : \JABG \longrightarrow \FLA
\EEN
of $\GAB$-schemes which sends  a closed point 
$\ZA$ of $\JABG$ to the partial flag $p_{\GA} \ZA$ determined by the filtration $\HT_{\bullet}$ at $\ZA$, i.e. we have
\BEN\label{p-ZA}
p_{\GA} \ZA]=[\scriptstyle{\HO{Z})=\HT_{-\LG-1} \ZA \supset \HT_{-\LG} \ZA \supset  \ldots \supset \HT_{-1} \ZA
 \supset \HT_0 \ZA=0}]\,,
\EEN
 for every $\ZA \in \JABG$.
\begin{defi}\label{period}
The morphism $p_{\GA}$ in (\ref{p}) is called the period map of $\JABG$.
\end{defi}

\begin{rem}\label{periods}
From the direct sum decomposition of $\HT_{-i}$ in (\ref{ordH-i}) it follows that the groups $\HH^p \ZA$, for $p=0,\ldots,\LG$, completely determine the value $p_{\GA} \ZA$, for
every $\ZA \in \JABG$. For this reason we often refer to the collection
 $\{\HH^p \ZA \}_{p=0,\ldots,\LG}$ as periods of $\ZA$.
\end{rem}

It turns out that the infinitesimal variation of periods along the fibres of the natural projection $\pi :\JAB_{\GA} \longrightarrow \GAB =\GA^{(0)}_{conf}$ in (\ref{pi-JAB})
 is closely related to the geometric properties of
the configurations $Z$ parametrized by $\GAB$. The point is that the relative differential of the period map $p_{\GA}$ can be computed entirely in terms of the multiplication morphisms
in (\ref{formalGr-H}).
\\
\\
{\bf 3.4. The relative differential of $p_{\GA}$.} Let ${\cal T}_{\pi}$ be the relative
tangent sheaf of the projection
$\pi :\JAB_{\GA} \longrightarrow \GAB = \GA^{(0)}_{conf}$ in (\ref{pi-JAB})  and let
${\cal T}_{Fl_{\GA}}$ be the relative tangent sheaf of the projection $Fl_{\GA} : \FLA \longrightarrow \GAB$.
If we let 
\begin{equation}\label{filt-U}
Fl^{\ast}_{\GA} (\FF_{\GAB}) = {\cal U}_{\LG+1} \supset {\cal U}_{\LG} \supset \cdots \supset {\cal U}_{1} \supset {\cal U}_{0} =0
\end{equation}
to be the universal partial flag of  $Fl^{\ast}_{\GA} (\FF_{\GAB})$ on $\FLA$, then the relative tangent sheaf ${\cal T}_{Fl_{\GA}}$ has the following description.
\begin{equation}\label{reltan-FL}
{\cal T}_{Fl_{\GA}} \cong  \bigoplus^{\LG+1}_{p=1} \HOM ({\cal U}_{p} /{\cal U} _{p-1}, Fl^{\ast}_{\GA} (\FF_{\GAB}) / {\cal U}_{p} )\,.
\end{equation}
Observe that by definition of the period map $p_{\GA}$ the pullback $p^{\ast}_{\GA} ({\cal U}_{p}) = \HT_{-p} \otimes \OO_{\JAB_{\GA}}$, for every $p =0,\ldots,\LG+1$.
For the rest of this discussion we will be working on $\JAB_{\GA}$, so all sheaves will be considered on this stratum. In particular, to simplify the notation, we will omit
tensoring with $\OO_{\JAB_{\GA}}$, whenever no ambiguity is likely.
   
 With the above preliminaries in mind, the relative differential $d_{\pi} ( p_{\GA})$ 
of $p_{\GA}$ can be written as follows 
$$
d_{\pi} ( p_{\GA}) : {\cal T}_{\pi} \longrightarrow p^{\ast}_{\GA} {\cal T}_{Fl_{\GA}} \cong \bigoplus^{\LG+1}_{p=1} \HOM (\HT_{-p} / \HT_{-p+1}, \FT / \HT_{-p} )\,.
$$ 
But in fact $d_{\pi} (p_{\GA})$ is subject to the Griffiths transversality condition and, moreover, we can calculate it using the multiplication in $\FT$.
 Before stating this result we need to discuss the relation between the relative tangent sheaf 
${\cal T}_{\pi}$ and the sheaf $\HT_{-1}$. 

This relation is not only central in the calculation of $d_{\pi} (p_{\GA})$, but it also becomes crucial
in the representation theoretic considerations of {\bf Part II}. The whole sense of this relationship can be summarized by saying that sections of
$\HT_{-1}$ can be viewed as {\it vertical}\footnote{`vertical' here and later on means vertical with respect to the morphism $\pi: \JABG \longrightarrow \GAB$.}\label{v}
vector fields on $\JABG$ and vice verse, vertical vector fields on $\JAB$ can be identified with sections of $\HT_{-1}$ and hence treated as functions on the
underling configurations of points in $X$.

The aforementioned relation between ${\cal T}_{\pi}$ and $\HT_{-1}$ is already implicit (on the level of closed points) in Lemma \ref{Htild}.
 Indeed, given a point $\ZA \in \JAB$, the fibre $\JAB_Z$ of $\pi: \JABG \longrightarrow \GAB$ over $[Z]$ is a non-empty Zariski open subset
of $\PP(\EZ)$. So we can identify the fibre 
${\cal T}_{\pi} \ZA$ of ${\cal T}_{\pi}$ at $\ZA$ with the quotient space $\EZ / \CC\{\alpha\}$.
On the other hand $\HT_{-1} \ZA =\HT \ZA$. This combined with the properties 2) and 3) of Lemma \ref{Htild} yield natural isomorphisms
$$
 \HT_{-1} \ZA / \CC \{1_Z \} \cong \EZ / \CC\{\alpha\} \cong  {\cal T}_{\pi} \ZA\,, 
$$
where $1_Z$ stands for the constant function of value $1$ on $Z$. The following is a sheaf version of the above isomorphism
between $\HT_{-1} \ZA / \CC \{1_Z \}$ and ${\cal T}_{\pi} \ZA $. 
\begin{pro}\label{H0-vertT}
Let $\OO_{\JAB_{\GA}} \subset \HT_{-1}$ be the natural inclusion from Proposition-Definition \ref{filtrn}, a).
Set $\HH := \HT_{-1} / \OO_{\JAB_{\GA}}$. Then one has a natural isomorphism
\begin{equation}\label{isoH0-vT}
M:\HH \stackrel{\cong}{\longrightarrow} {\cal T}_{\pi}.
\end{equation} 
\end{pro} 
The identification in (\ref{isoH0-vT}) gives another perspective on the sheaf $\HT (=\HT_{-1})$.
Namely, the obvious exact sequence for $\HT$ , the first row below, 
\BEN\label{exseq-HT}
\xymatrix{
0\ar[r]  & \OO_{\JAB_{\GA}} \ar@{=}[d] \ar[r]& \HT \ar[d]^{\tilde{M}}_{\cong} \ar[r]& \HH \ar[d]^{M}_{\cong} \ar[r]& 0 \\
0\ar[r]  & \OO_{\JAB_{\GA}} \ar[r]  & {\mathcal{D}}^{\scriptscriptstyle{\leq1}}_{\pi} \ar[r] &{\cal T}_{\pi}\ar[r]& 0 }
\EEN
gets identified with the sequence for the sheaf ${\mathcal{D}}^{\scriptscriptstyle{\leq1}}_{\pi}$ of the {\it relative} first order differential operators
on $\JABG$, the second row in (\ref{exseq-HT}). The isomorphism $\tilde{M}$ in the middle of (\ref{exseq-HT}) follows from the fact that  the top
exact sequence in (\ref{exseq-HT}) admits a canonical splitting. 
\begin{pro}\label{HT-split}
Let $\HH^{\prime} = \big{(}\OO_{\JABG} {)}^{\perp}  \cap \HT$ be the orthogonal complement in $\HT$ of the subsheaf $\OO_{\JABG} \subset \HT$  with respect to the quadratic form $\QT$ in (\ref{q}).
Then one has the orthogonal decomposition
\BEN\label{f:splitHT}
\HT =\OO_{\JABG} \oplus \HH^{\prime}\,.
\EEN  
Furthermore, the epimorphism in the top exact sequence in (\ref{exseq-HT}) induces an isomorphism
\BEN\label{Mpr}
M^{\prime} :\HH^{\prime} \stackrel{\cong}{\longrightarrow}  {\cal T}_{\pi}\,.
\EEN
\end{pro}
 The isomorphism $\tilde{M}$ in (\ref{exseq-HT}) allows to think of sections of $\HT$ as  relative first order differential operators on $\JABG$.
This indicates a possibility that the sheaf of relative differential operators
${\mathcal{D}}_{\pi}$ of $\JABG$ could provide a non-commutative deformation of the direct sum decomposition (\ref{ordFT}) and lead eventually to {\it non-commutative}
Hilbert schemes of points of projective surfaces. This is something for future to tell.

We now return to the relative differential $d_{\pi} (p_{\GA})$ and state its main properties.
\begin{pro}\label{pro-Grif-trans}
The relative differential $d_{\pi} (p_{\GA})$ satisfies Griffiths transversality condition
\begin{equation}\label{Grif-trans}
d_{\pi} (p_{\GA}) : {\cal T}_{\pi} \longrightarrow  
\bigoplus^{\LG}_{m=1} \HOM (\HT_{-m} /{ \HT_{-m+1}}, \HT_{-m-1} /{ \HT_{-m} }) 
\end{equation}                                                                                                                                                            
Furthermore, let
 $d_{\pi} (p_{\GA})_m$ be the $m$-th component of 
 $d_{\pi} (p_{\GA})$. Then for any local section $v$ of $ {\cal T}_{\pi}$ and any local section
$h$ of $\HT_{-m}$ one has the following:
\BEN\label{diff-formula}
d_{\pi} (p_{\GA})_m (v) (h) \equiv -m  \tilde{v} \cdot h\, (mod\, \HT_{-m})
\EEN
where $\tilde{v}$ stands for an arbitrary lifting of $M^{-1} (v)$, where $M$ is the isomorphism in (\ref{isoH0-vT}),
 to a local section of $\HT$
and $\tilde{v} \cdot h $
 stands for the product in $\FT$ and it is independent of a lifting chosen after factoring out by $ \HT_{-m}$.
\end{pro}

The first immediate consequence of the formula (\ref{diff-formula}) is that the subspace $\HT_{-\LG} \ZA$, the fibre of the subsheaf $\HT_{-\LG}$ in (\ref{filtHT-JG}), is a subring\footnote{
the subring  $\HT_{-\LG} \ZA$ can be described using the morphism $\kappa_{\ZA}$ in (\ref{kappaZ}): let $Z^{\prime}$ be the image of $\kappa_{\ZA}$, then 
$$\HT_{-\LG} \ZA =\kappa^{\ast}_{\ZA} (H^0(\OO_{Z^{\prime}}))$$
is the pullback of the space of functions on $Z^{\prime}$.} of $\HO Z)$, and it
is independent of $[\alpha]$ in the fibre $\JAB_Z$ of $\JABG$ over $[Z]$.
In particular, this implies that the sheaf $\HT_{-\LG}$ is the pullback of a sheaf on $\GAB$ and that the interesting part of the period map $p_{\GA}$ concerns only the decomposition
\BEN\label{ordLG}
\HT_{-\LG} =\bigoplus^{\LG -1}_{p=0} \HH^p .
\EEN
  
 From the conceptual point of view, the formula (\ref{diff-formula}) expressing the (relative) differential in terms of the multiplication in the sheaf of rings
$\FT$ is useful, because it translates the diffeo-geometric properties of the period map $p_{\GA}$ into algebro-geometric properties of 
the configurations parametrized by $\GAB$. For example, the failure of $p_{\GA}$ to be injective on the fibre
$\JAB_Z$ of $\JAB_{\GA}$ over a point $[Z] \in \GAB$ (i.e. the failure of Torelli property for $p_{\GA}$)  can be interpreted as a {\it canonical}
decomposition of $Z$ into disjoint union of subconfigurations.\label{p-Zdec} If the injectivity of $p_{\GA}$ fails over a general point of $\GAB$, then one has the decomposition
for every $Z$, with $[Z] \in \GAB$, and this leads to a non-trivial factorization of the covering map
$p_{2}: \ZD_{\GAB} \longrightarrow \GAB$, where $p_2$ is as in (\ref{uc}) and $\ZD_{\GAB} = p^{-1}_2 (\GAB)$ is the part of the universal scheme $\ZD$ lying over $\GAB$.
Thus it should be clear that the period map in (\ref{p}) is useful for revealing geometric properties of $0$-dimensional subschemes of $X$.
   
 The utility of the decomposition (\ref{ordFT}) for defining the period map $p_{\GA}$ turns out to be only a part of the story. 
The same decomposition together with the multiplicative structure of 
$\FT$ are crucial in defining
the sheaf of Lie algebras {\boldmath$\GT_{\mbox{\unboldmath${\GA}$}}$}. 
\\
\\
{\bf 3.5. The sheaf of the Lie algebras {\boldmath$\GT_{\mbox{\unboldmath${\GA}$}}$.} }\label{Lie}
To attach Lie algebras to points of $\JA$ we view local sections of the sheaf
$\HT =\HT_{-1}$
in (\ref{H-filt}) as operators of the multiplication in  the sheaf of rings $\FT$, i.e. we consider the inclusion
\begin{equation}\label{m-D}
\xymatrix@1{
{D: \HT} \ar[r]& {\ENDO(\FT)}  }
\end{equation}
which sends a local section $t$ of $\HT$ to the operator
$D(t)$ of the multiplication by $t$ in $\FT$.
\\
\indent
Over the configuration subset
$\GAB=\GAC$ of an admissible component $\GA \in \CSA$
we have defined the subscheme
$\JAB_{\GA}$ (see (\ref{pi-JAB}) for notation),
where the orthogonal decomposition
(\ref{ordFT}) holds. Using this decomposition we write the multiplication operator
$D(t)$ in the ``block" form, where the block 
$D_{pm} (t) $ is the local section of
$\HOM (\HH^p, \HH^m)$ obtained by taking the restriction
$D_p (t)$ of $D(t)$ to the summand $\HH^p$ of the decomposition in (\ref{ordFT}) followed by the projection onto its $m$-th summand.

It turns out that for every $p$ one has at most three non-zero blocks
\begin{equation}\label{d-Dpt}
D_p(t) = D^{-}_{p} (t) + D^0_{p} (t)  + D^{+}_{p} (t)
\end{equation}
where 
$D^{\pm}_p (t)$ (resp. $D^0_{p} (t)$) denote the blocks $D_{p,p\pm1}$ (resp. $D_{p,p} (t)$).
This implies that the multiplication operator $D(t)$ admits the following decomposition
\begin{equation}\label{d-Dt}
D(t) = D^{-} (t) + D^0 (t)  + D^{+} (t)
\end{equation}
where
$\displaystyle{D^{\pm} (t) =\sum^{\LG}_{p=0} D^{\pm}_p (t)}$ and $\displaystyle{D^0 (t) =\sum^{\LG}_{p=0} D^{0}_p (t)}$.
In particular, the operators $D^{\pm} (t)$ are of degrees $\pm 1$, while 
$ D^{0}(t) $
is of degree $0$, with respect to the grading in (\ref{ordFT}).
Thus on $\JAB_{\GA}$ the morphism $D$ in (\ref{m-D}) admits the triangular decomposition
\begin{equation}\label{d-D}
 D = D^{-}  + D^{0}  + D^{+}
\end{equation}
and we define 
{\boldmath$\GT_{\mbox{\unboldmath${\GA}$}}$}
to be the subsheaf of Lie subalgebras of
$\ENDO(\FT)$
generated by the subsheaves
$ D^{\pm} (\HT)$ and $D^{0} (\HT)$.
\\
\\
\indent
The sheaf of Lie algebras {\boldmath$\GT_{\mbox{\unboldmath${\GA}$}}$}
could be viewed as a nonabelian counterpart of the (abelian) Lie structure of the classical Jacobian.
But there is more to it than a simple analogy since 
{\boldmath$\GT_{\mbox{\unboldmath${\GA}$}}$}
allows to bring the powerful methods of the representation theory of Lie algebras and Lie groups to address various geometric questions.
This is the contents of \cite{[R2]} and we will return to this in more details in the second part of this survey. For now we will pursue
another aspect of the triangular decomposition in (\ref{d-D}) related to the ideas of Simpson in \cite{[S]} on nonabelian Hodge theory.
Namely, we interpret the morphism $D$ in (\ref{m-D}) as a Higgs structure on $\FT$ and then use the equation (\ref{d-D}) to produce
a distinguished family of such structures.
\section{Relative Higgs structures on $\FT$ and the associated Fano varieties}
 
{\bf 4.1. Higgs structures on $\FT$.} We begin by rewriting the morphism $D$ in $(\ref{m-D})$ as follows
 \begin{equation}\label{mult1}
D:\FT  \longrightarrow  \HT^{\ast} \otimes \FT   
 \end{equation}
 To fix our notation and terminology we will need some
 generalities about morphism written in this form.
 \\
 \\
 \indent
 Let ${\mathcal{M}}$ and ${\mathcal{N}}$ be two vector bundles over a scheme
 $S$ (as usual we make identification of vector bundles over $S$
 and locally
 free $\OO_S$-modules) and let
 \begin{equation}\label{mor}
A:{\cal M} \longrightarrow   {\cal N}^{\ast} \otimes {\cal M}   
 \end{equation}
 be a morphism of $\OO_S$-modules. For a section $n$ of 
 ${\cal N}$
 over an open subset $U \subset S$ we denote
 $$
A(n):{\cal M}\otimes \OO_U \longrightarrow   {\cal M}\otimes \OO_U   
$$
the endomorphism of ${\cal M}\otimes \OO_U$ induced by $A$
(we often omit the reference to an open subset in the above notation).

Given two morphisms $A,\, B:{\cal M} \longrightarrow {\cal N}^{\ast} \otimes {\cal M}$
 we define
$$
A\wedge B:{\cal M} \longrightarrow   \wedge^2{\cal N}^{\ast} \otimes {\cal M}   
$$
as follows:
\begin{equation}\label{AB-def}
(A\otimes id_{{\cal N}^{\ast}}) \circ B - (B\otimes id_{{\cal N}^{\ast}} )\circ A 
\end{equation}
\begin{rem}\label{AB-formula}
\begin{enumerate}
\item[1).]
We write $A^2$ for $A\wedge A$.
 \item[2).] 
For a local section $n \wedge n^{\prime}$ of $\wedge^2 {\cal N}$
the morphism $A \wedge B$ is given by the following formula
$$
(A\wedge B )(n \wedge n^{\prime}) = [A(n^{\prime}),\,\, B(n)]
$$
where the bracket is the commutator of endomorphisms. In particular,
$A^2 = 0$ if and only if $A(n)$ and $A(n^{\prime})$ are commuting endomorphisms
for any local sections $n ,\, n^{\prime}$ of ${\cal N}$.
\end{enumerate}
\end{rem}
\begin{defi}\label{Higgs-mor}
Let $A$ be as in $(\ref{mor})$. It is said to be a Higgs endomorphism of 
${\cal M}$
with values in ${\cal N}^{\ast}$ if $A^2= 0$.
\end{defi}
\begin{rem}\label{Higgs}
In our terminology a Higgs bundle over $S$ (see \cite{[S]}) is a bundle with a 
Higgs endomorphism having its values in the cotangent bundle of $S$.
More generally, let $f: S \longrightarrow B$ be a smooth morphism of relative 
dimension
$\geq 1$. We say that a bundle
${\cal M}$ over $S$ is a relative Higgs bundle if it has a Higgs endomorphism
with values in the relative cotangent bundle of $f$. In this case a Higgs endomorphism
of ${\cal M}$ will be called a relative Higgs field.
\end{rem}
We will now return to our considerations of the multiplicative action of $\HT$ on $\FT$.
Since the multiplication in $\FT$ is commutative it follows immediately 
\begin{lem}\label{D-Higgs}
The morphism
$D$
in (\ref{mult1}) is a Higgs morphism of $\FT$ with values in $\HT^{\ast}$.
In particular, $D^2 =0$ and the decomposition in (\ref{d-D}) is subject to the following identities
\begin{enumerate}
\item[(i)]
$D^2 = (D^{-})^2= (D^{+})^2 = 0$,
\item[(ii)]
$D^{-} \wedge D^{0} + D^{0}\wedge D^{-} =
D^{+} \wedge D^{0} + D^{0}\wedge D^{+} = 0$,
\item[(iii)]
$(D^{0})^2 + D^{-} \wedge D^{+} + D^{+} \wedge D^{-} = 0$.
\end{enumerate}
\end{lem}
\begin{rem}\label{relHiggs}
The orthogonal splitting of $\HT =\OO_{\JABG} \oplus \HH^{\prime}$ in (\ref{f:splitHT}) implies that a Higgs structure with values in $\HT^{\ast}$ is automatically
a Higgs structure with values in $\HH^{\prime \ast}$. This and the identification $M^{\prime} :  \HH^{\prime} \cong {\cal T}_{\pi}$ in (\ref{Mpr}) imply that every Higgs structure with values in $\HT^{\ast}$
determines a relative Higgs structure, i.e. the one with values in the relative cotangent bundle ${\cal T}^{\ast}_{\pi}$ of $\pi: \JABG \longrightarrow \GAB$.
\end{rem}  

  A simple consequence of Lemma $\ref{D-Higgs}$ is that the components
$D^0_p$, $D^{\pm}_p$ of $D^{0}$ and $D^{\pm}$ satisfy the following relations
\begin{eqnarray}\label{relp}
D^{\pm}_{p\pm1} \wedge D^{\pm}_p =0\,,   \nonumber\\
D^{\pm}_p \wedge D^0_p + D^0_{p\pm 1} \wedge D^{\pm}_p =0\,,\\
(D^0_p)^2 +  D^{-}_{p+1} \wedge D^{+}_p + D^{+}_{p-1} \wedge D^{-}_p = 0\,, 
\nonumber
\end{eqnarray}
for every $p= 0,\ldots,\LG-1$ and where we use the convention that
$D^{\pm}_p $ (resp. $ D^0_p$) is zero whenever the index $p$ is not in the
above range. Using these relations we can create a large family of Higgs structures on $\FT$.

This is achieved
by taking sufficiently general deformation
of $D$. The construction we are about to discuss requires the weight $\LG$ of the orthogonal
decomposition of $\FT$ to be $\geq 2$. The case $\LG=1$ is very special. It implies in particular, that the components $D^{\pm}$
of the decomposition (\ref{d-D}) are both equal to zero. Thus in this case the sheaf $\LAGT$ is a sheaf of abelian Lie algebras
and the configurations $Z$ parametrized by such a component $\GA$ have a very special geometry. The phenomenon here is somewhat reminiscent of
the hyperellipticity in the theory of curves and will be considered elsewhere. 
So from now on we assume $\LG \geq2$.

 Another remark which is in order is that the operators of the multiplication $D(t)$, for any local section $t$ of $\HT$, preserve the subsheaf $\HT_{-\LG}$ of the filtration
 in (\ref{filtHT-JG}).\footnote{this follows from the definition of
$\HT_{-\LG}$.} This implies that $D^{+}_{\LG-1} (t) =0$ and in what follows we consider the restrictions of all operators to the subsheaf 
\BEN\label{ordLG1}
\HT_{-\LG} = \bigoplus^{\LG-1}_{p=0} \HH^p
\EEN
 of $\FT$. 
\\
\\
{\bf 4.2. Nonabelian $(1,0)$-Dolbeault variety.}
 We consider a sufficiently general
deformation of $D$ of the following form
$$
\sigma(t,x,y) = \sigma^0(t) +  \sigma^{+} (x) + \sigma^{-} (y)
$$
where
\begin{equation}\label{txy}
\sigma^0(t) = \sum^{\LG-1}_{p=0} t_p D^0_p ,\,\,\,\,\,
\sigma^{+} (x) = \sum^{\LG-2}_{p=0} x_p D^{+}_p ,\,\,\,\,\,
\sigma^{-} (y) = \sum^{\LG-2}_{p=0} y_{p} D^{-}_{p+1}
\end{equation}
and
$t=( t_p) \in {\bf C^{\LG}},\,\,x=(x_p),\,\,y=(y_{p}) \in {\bf C^{\LG-1}} $ are deformation parameters.
One derives sufficient conditions for the morphism $\sigma(t,x,y)$
to be Higgs by writing out the expansion
\begin{eqnarray*}
\sigma^2(t,x,y) = \sigma^0 (t) \wedge \sigma^{+} (x)
 + \sigma^{+} (x)\wedge\sigma^0(t)+(\sigma^0(t))^2 +
\sigma^{+} (x)\wedge \sigma^{-} (y)&\\
 + \sigma^{-} (y)\wedge\sigma^{+} (x) + \sigma^0(t) \wedge \sigma^{-} (x) +
\sigma^{-} (x)\wedge\sigma^0(t)&
\end{eqnarray*}
according to the degree with respect to the grading in (\ref{ordLG1}) of components of $\sigma(t,x,y)$ . Then 
$\sigma^2(t,x,y) =0$ yields the vanishing in each degree
\begin{equation}\label{sys}
\left\{
\begin{array}{ll}
\sigma^0(t) \wedge \sigma^{+} (x) + 
\sigma^{+} (x)\wedge\sigma^0(t)& =0 \\ 
(\sigma^0(t))^2 + \sigma^{+} (x)\wedge\sigma^{-} (y) + 
\sigma^{-} (y)\wedge\sigma^{+} (x)& =0\\
\sigma^0(t) \wedge \sigma^{-} (x) +
\sigma^{-} (x)\wedge\sigma^0(t)& =0
\end{array}
\right.
\end{equation}
Substituting in the expressions of (\ref{txy}) and using the relations
(\ref{relp}) we arrive at the following system of equations
\begin{equation}\label{sys1}
\left\{
\begin{array}{lll}
x_p (t_{p+1} - t_p)& =&0\\
y_p (t_{p+1} - t_p)& =&0\\
x_p y_{p} - t^2_p& = & x_p y_{p} - t^2_{p+1} =0
\end{array}
\right.
\end{equation}
for $p= 0,\ldots,l-2$.
This yields the set of solutions
\begin{equation}\label{sol}
\hat{H} = \left\{ \left.(z,x,y)\in {\bf C\times C^{\LG-1}\times C^{\LG-1}} \right|
x=(x_p),\,\,y=(y_p),\,\, x_p y_p =z^2,\,\,p=0,\ldots,\LG-2 \right\}\,.
\end{equation}
For $(z,x,y)\in \hat{H}$ denote by 
$\displaystyle{
\sigma(z,x,y) = zD^0 + \sum^{\LG-2}_{p=0} x_p D^{+}_p + 
\sum^{\LG-2}_{p=0} y_p D^{-}_{p+1}}$
the corresponding Higgs morphism. It is clear that scaling 
$\sigma(z,x,y)$ by $\lambda \in{\bf C^{\ast}}$ gives a $\bf C^{\ast}$-
action on $\hat{H}$. Furthermore, conjugating
$\sigma(z,x,y)$ by an automorphism 
$$
g= \sum^{\LG-1}_{p=0} g_p id_{\bf H^p},\,\, 
(g_p\in{\bf C^{\ast}},\,p=0,\ldots,\LG-1)
$$
of $\FT$ gives a gauge equivalent Higgs morphism 
\begin{equation}\label{conj}
g \sigma(z,x,y) g^{-1} = zD^0 + \sum^{\LG-2}_{p=0}\frac{g_{p+1}}{g_p} x_p D^{+}_p
 + 
\sum^{\LG-2}_{p=0}\frac{g_{p}}{g_{p+1}} y_p D^{-}_{p+1}\,.
\end{equation}
All together this defines an action of the torus
$\hat{S}= {\bf (C^{\ast })^{\LG}}$ on the variety $\hat{H}$. From
(\ref{conj}) we deduce that the action has the
following form. For 
$\tau = (\lambda,\lambda_0,\ldots,\lambda_{\LG-2}) \in \hat{S}$ and
$(z,x,y) \in \hat{H}$ we have
\begin{equation}\label{T-act}
\tau \cdot (z,x,y) = \lambda (z, \lambda_0 x_0,\ldots,\lambda_{\LG-2} x_{\LG-2},
 \lambda^{-1}_0 y_0,\ldots, \lambda^{-1}_{\LG-2} y_{\LG-2} )\,. 
\end{equation}
If we factor out $\hat{H}$ by the scaling $\bf C^{\ast}$-action we obtain
the projectivization of $\hat{H}$ which will be denoted by $H$.
\begin{defi}\label{Dol}
$H$ is a variety of the homothety equivalent non-zero Higgs morphisms of
$\FT$. This variety will be called nonabelian $(1,0)$-Dolbeault\footnote{in \RI, \S4, this variety was called `nonabelian Albanese'. This terminology is not quite appropriate, since classically, the
Albanese variety involves taking the {\it dual} of the space $H^{1,0}$ of holomorphic $1$-forms.
 The variety $H$ constructed here is certainly more like a direct analogue of the space of holomorphic $1$-forms itself, if one thinks of Higgs morphisms as a nonabelian version of holomorphic $1$-forms. Hence the change of terminology.} variety associated to ${\JAB}_{\GA}$.
\end{defi}
From (\ref{sol}) one easily obtains the following projective description of $H$.
\begin{pro}\label{H1}
Let $\PP^{2(\LG-1)}$ be a projective space with the coordinates
$T,X_p,Y_p$, $(p=0,\ldots,\LG-2)$. Then $H$ is a complete intersection of $\LG-1$
quadrics $X_p Y_p =T^2 \,(p=0,\ldots,\LG-2)$ in $ \PP^{2(\LG-1)}$.
In particular, $H$ is a Fano variety of dimension $(\LG-1)$ and degree $2^{\LG-1}$ with the dualizing
sheaf $\omega_H = \OO_H (-1)$.
\end{pro}

From the adjunction formula one deduces the following.
\begin{cor}\label{h-sec}
The hyperplane sections of $H$ are Calabi-Yau varieties of dimension
$\LG-2$.
\end{cor}

The  $(1,0)$-Dolbeault variety $H$ comes together with the distinguished divisor
$H_0$ corresponding to the hyperplane section defined by $T=0$.
This divisor corresponds to Higgs morphisms having components
 of type $D^{\pm}_p $ only. Projectively, the divisor $H_0$ is a degenerate divisor - it is 
a union of projective spaces.
To give its precise projective description 
 we begin with
a more
intrinsic definition of the projective space in Proposition \ref{H1}.

Introduce the symbols\footnote{one could think of these symbols as 
labels of vertices of a certain graph - this is the trivalent graph we alluded to in the Introduction and which will be discussed in \S5.}
$V_0, V^{\pm}_p,\,p=0,\ldots,\LG-2$,  which will be eventually colored
by $D^0$ and $D^{\pm}_p,\,p=0,\ldots,\LG-2$, respectively, and let
$V$ be $\bf C$-vector space generated by these symbols. Let
$V^{\ast}$ to be the space dual to $V$ and define
$T,X_p,Y_p,\,(p=0,\ldots,\LG-2)$ to be its dual basis. Now the projective
space $\PP^{2(\LG-1)}$ in Proposition \ref{H1} is just $\PP(V)$ and the divisor 
$H_0$ is defined by the equations $T=0$ and $X_p Y_p =0 \,(p=0,\ldots,\LG-2)$.
Then we have the following description of $H_0$.
\begin{lem}\label{H0}
For any subset $A$ of the set of indicies $I=\{0,\ldots,\LG-2\}$ let
$\hat{\Pi}_A$ be  the subspace of $V$ spanned by
the vectors $\{V^{+}_i,V^{-}_j \mid i\in A,j\in I\setminus A \}$ and let $\Pi_A$
 be its projectivization.
Then $\displaystyle{H_0 = \bigcup_A \Pi_A}$, where the union is taken over
all subsets $A$ of the set $I$.
\end{lem}
From the definition of the irreducible components $\Pi_A$ of $H_0$ and
the equations defining $H$ (Proposition \ref{H1}) 
 one easily obtains the following.
\begin{lem}\label{singH}
\begin{enumerate}
\item[(i)]
 For any two subsets $A,B$ of $I$ the intersection
$ \Pi_A \cap \Pi_B$ is the projectivization of the vector space
spanned by the set 
$\{V^{+}_i,V^{-}_j \mid i\in A\cap B,\,\, j\in A^-\cap B^- \}$,
where $A^-$ denotes the complement of $A$ in $I$.
\item[(ii)]
The singularity locus of $H$ is 
$ \displaystyle{Sing(H) = \bigcup_{A\not=B} \Pi_A \cap \Pi_B }$,
where the union is taken over all pairs  of distinct subsets $A,B$ of the index
set $I$.
\end{enumerate}
\end{lem}
We have seen in Corollary \ref{h-sec} that the hyperplane sections of
$H$ are Calabi-Yau varieties. One can argue that the divisor
$H_0$ is degenerate from this perspective as well. Namely,
the divisor $H_0$ comes with a {\it degenerate} symplectic structure - the projective spaces composing it are
the projectivized  Lagrangian subspaces as explained in the following
\begin{lem}\label{Lag}
Let $V^{\prime}$ be the subspace of $V$
spanned by the vectors $V^{\pm}_p,\,(p=0,\ldots,\LG-2)$. It admits a
natural symplectic structure $\omega$ with respect to which the subspaces
$\hat{\Pi}_A$ in Lemma \ref{H0} are Lagrangian subspaces of the symplectic
space $(V^{\prime},\omega)$. 
\end{lem}

\begin{defi}\label{Lag1}
The divisor $H_0$ in Lemma \ref{H0} is called the Lagrangian cycle
of $H$ and its irreducible components ${\Pi}_A$ are called
 Lagrangian manifolds of $H$.
\end{defi}

From (\ref{T-act}) it follows that the nonabelian $(1,0)$-Dolbeault variety $H$ admits a natural torus action which gives it a structure of a
toric
variety.
\begin{pro}\label{H2-0}
The nonabelian $(1,0)$-Dolbeault variety $H$ is a projective toric variety
with an action of the torus $S = \bf (C^{\ast})^{\LG-1}$. Its
fan $\Delta$ is the fan in $\bf R^{\LG-1}$ generated by the vertices of the
cube $[-1,1]^{\LG-1}$ and the vertices of the cube are in
1-to-1 correspondence with the irreducible components of the divisor
$H_0$. In particular, $Pic(H)$ is generated by the irreducible
components of the ``Lagrangian" cycle $H_0$.
\end{pro}

It should be observed that 
the nonabelian $(1,0)$-Dolbeault variety $H$ and its Lagrangian subcycle $H_0$ depend only on the weight
$\LG$ of the direct sum decomposition in (\ref{ordFT}) and the decompositions in (\ref{d-Dpt}) that was used to define
 the Higgs morphisms parametrized by $H$. These data can be conveniently represented and determined by a certain trivalent
graph which will be discussed next.
\section{The trivalent graph associated to $\JAB_{\GA}$}
{\bf 5.1. Definition of the graph.} The graph $G=G(\JAB_{\GA})$ which we associate to $\JAB_{\GA}$ is basically a pictorial
representation of the orthogonal decomposition 
(\ref{ordFT}) and the action of the morphisms
$D^0,D^{\pm}$ in (\ref{d-D}) on the summands of this decomposition.
More precisely, we take  $\LG$ parallel vertical edges with
upper (resp. lower) vertices aligned on a horizontal line (see (\ref{G})). These
vertical segments should be thought of as a pictorial representation of 
the morphism $D^0$ preserving the summands of the decomposition (\ref{ordFT}). The segments are naturally
ordered, from left to right, by the index set
$I=\{0,1,\ldots,\LG-1\}$. The vertices of the $i$-th edge will be labeled by
$i^u$ and $i^d$, for the upper and lower one respectively, and they should be viewed as a pictorial representation of the $\LG$ summands
 of the direct sum decomposition in (\ref{ordLG}).\footnote{in the paragraph just above (\ref{ordLG}) it is explained that the interesting part of the period map $p_{\GA}$
is contained in $\HT_{-\LG} =\displaystyle{\bigoplus^{\LG -1}_{p=0} \HH^p}$. This is why we are ignoring the last summand $\HH^{\LG}$ in (\ref{ordFT}).}

 We now connect
the neighboring vertices as follows. For every $i\in I$ connect 
$i^u$ to $(i-1)^d $ and $(i+1)^d$
using the convention that $(-1)^d = (\LG-1)^d$ and $l^d_{\GA} = 0^d$.
This gives the following trivalent graph
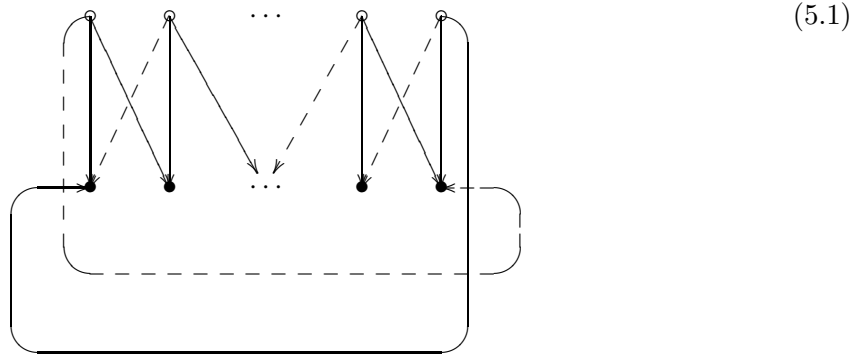
\begin{equation}\label{G}
\xymatrix{
&*=0{\circ} \ar[2,0] \ar[2,1] \ar@{--{>}}`l[d]                               
                               `[3,1]                               
                               `[2,5]                               
                                `[2,4]
                                [2,4]
                              &
                               *=0{\circ}\ar[2,0] \ar[rdd] \ar@{--{>}}[ddl]&
\cdots&*=0{\circ} \ar[2,0] \ar@{--{>}}[ldd] \ar[2,1]&
*=0{\circ} \ar[2,0] \ar@{--{>}}[2,-1] \ar`r[d] 
                              `[4,-2]
                              `[2,-5]                              
                              `[2,-4]
                              [2,-4]&\\
& & & & & & \\
&*=0{\bullet}&*=0{ \bullet}&\cdots&*=0{\bullet}&*=0{\bullet}&\\
& & & & & &\\
& & & & & &}
\end{equation}
which will be denoted $G(\JAB_{\GA})$ or simply by $G$ if no ambiguity is likely.
We agree that all edges of $G$ are oriented from ``up" to ``down".
The graph $G$ with this orientation will be denoted by $\vec{G}$. In this 
orientation the edges incident to an upper (resp. lower) vertex of $G$ are
always out (resp. in)-going. At every upper vertex of $G$ we fix a 
counterclockwise
order of incident edges and color them, starting with the vertical one,
by ``$0$", ``$+$", and ``$-$", respectively. This will be called the 
{\it natural coloring of } $G$. It is chosen to correspond to the 
morphisms $D^0, D^{+}, D^{-}$ and the orientation of edges of $\vec{G}$
matches the sense of action of these operators on a vector of pure
degree $p$ (in the decomposition (\ref{ordLG})) placed at the vertex $p^u$ of $G$.

At this stage it should be clear that the graph $G(\JAB_{\GA})$ captures the form, i.e. the number of summands,  of the direct sum decomposition
in (\ref{ordLG}) as well as the action of operators generating the sheaf of Lie algebras {\BM$\LAGT$}. If we want, in addition, to record the ranks of the summands in (\ref{ordLG}), then
we assign weights to the vertices of $G(\JAB_{\GA})$ by the rule
\BEN\label{weight-f}
(p^u,p^d) \mapsto h^p_{\GA} =rk (\HH^p)\,,
\EEN
for every $p$ in $\{0,\ldots,\LG-1\}$. In other words we view the {\it reduced Hilbert vector}
$h^{\prime}_{\GA} = ( h^0_{\GA},\ldots, h^{\LG-1}_{\GA})$ as a weight function of $G(\JAB_{\GA})$ and the pair
$(G(\JAB_{\GA}),h^{\prime}_{\GA}) $ as a {\it weighted graph}.   
\\
\\
{\bf 5.2. The graph $G(\JAB_{\GA})$ and quantum properties of $\JA$.}
Heuristically, the Jacobian $\JA$ is a ``quantum" object in the sense that it reveals very involved inner structure of points (particles) on
projective surfaces. From this perspective the graph $G(\JAB_{\GA})$ is a precise mechanism which transforms classical observables = functions 
on a configuration $Z$, with $\ZA \in \JAB_{\GA}$, to quantum observables = operator-valued generating series. Let us briefly explain this point.

Let $\ZA$ be a point of $\JAB_{\GA}$ and let $f \in \HO Z )$ be a function on $Z$. The direct sum (\ref{ordFT}) at the point $\ZA$ yields the following decomposition 
\begin{equation}\label{d-f}
f= \sum^{\LG}_{p=0} f_p  
\end{equation}
where $f_p \in \HH^p \ZA$ is the projection of $f$ on the $p$-th summand of the direct sum (\ref{ordFT}) at $\ZA$. We take the component
$f_0 \in \HH^0 \ZA = \HT \ZA$ and exploit the triangular decomposition in (\ref{d-Dt}) to obtain
\begin{equation}\label{d-Df0}
D(f_0)  = D^{-} (f_0) + D^0 (f_0)  + D^{+} (f_0)
\end{equation}
Now we insert the vector $f_0$ into the upper vertices of the graph $G(\JAB_{\GA})$ and replace the coloring ``$0$", ``$+$", ``$-$" of edges
of the graph by the operators $D^0 (f_0)$, $D^{+} (f_0)$, $D^{-} (f_0)$ respectively. This procedure colors all edges of 
$G(\JAB_{\GA})$ by elements of the Lie algebra
${\bf \tilde{\goth g}}\ZA$, the fibre of {\boldmath$\GT_{\mbox{\unboldmath${\GA}$}}$}
at $\ZA$. As a consequence we can assign to any path $\gamma$ of the graph $G(\JAB_{\GA})$ the operator which is the product of operators coloring the edges composing 
the path $\gamma$. All these path operators can be organized in a generating series. For more details and applications of this path operator technique see \RI, \S6. 
\\
\\
{\bf 5.3. From $G(\JAB_{\GA})$ to the Dolbeault variety $H$.} The graph $G(\JAB_{\GA})$ also allows to recover the nonabelian $(1,0)$-Dolbeault variety $H$ and its Lagrangian cycle in a purely formal way. Indeed,
let us label the upper vertices $i^u$, for $i=0,\ldots,\LG-2$, by $V^+_i$. Then we move down from each $i^u$ along the edge of $G(\JAB_{\GA})$ colored by `$+$'
to arrive
to the vertex $(i+1)^d$ and assign to it the label $V^{-}_i$. The set of labels $V^{\pm}_i$ $(i=0,\ldots,\LG-2)$ is completed by the ``edge" label $V_0$.
We consider the labels $V_0, V^{\pm}_0,\ldots, V^{\pm}_{\LG-2}$ as a basis of $\CC$-vector space which will be denoted by $V$. The projective space
$\PP(V)$ of dimension $2(\LG-1)$ is the one appearing in Proposition \ref{H1}. The homogeneous coordinates $T,X_p,Y_p, (p=0,\ldots,\LG-2)$ of that proposition
can be now intrinsically defined as the basis of $V^{\ast}$ dual to the basis $V_0, V^{\pm}_0,\ldots, V^{\pm}_{\LG-2}$ of $V$.
The variety $H$ is now (re)defined as the complete intersection of $\LG-1$ quadrics
$$
X_p Y_p = T^2, \,\,for\,\,p=0,\ldots,\LG-2\,,
$$
while the Lagrangian cycle $H_0$ is the divisor of $H$ defined by the hyperplane $T=0$. The discussion about $H_0$ in \S4, once the graph
$G(\JAB_{\GA})$ is in place, assumes  its full meaning.
 
It should be also clear now that the variety $H$ completely ``forgets" our surface $X$ and the Jacobian $\JA$.
In such a situation one expects to have correspondences between $\JAB_{\GA}$ and $H$.
\\
\\
{\bf 5.4. Geometric correspondence between $\JAB_{\GA}$ and $H$.}
Let us begin with a general remark about the nature of correspondences between 
$\JAB_{\GA}$ and the Dolbeault variety $H$. There are two types of correspondences we have in mind.
The first one is in line with the classical correspondence
between the Jacobian of a curve and its Albanese given by the Albanese map.
This map essentially involves integrating holomorphic $1$-forms on the Jacobian along a path.
The geometric correspondence we will be discussing here is of this nature.

The second type of correspondences is of completely different nature. It is based on the general philosophy of 
the Homological Mirror Symmetry conjecture of Kontsevich, \cite{[K]}. The properties of $H$ and its Lagragian cycle $H_0$ 
discussed in \S4, suggest that the pair $(H,H_0)$ should play a role of the ``symplectic" side of the mirror duality and should carry some kind of
Fukaya-type category, while
the pair $(\JAA_{\GA}, \TE_{\GA})$ should be its algebraic/holomorphic side and it should be equipped with a suitable subcategory of coherent sheaves.
 Then correspondences between $\JAA_{\GA}$ and $H$ should be functors between these categories. We will not pursue this direction here any further 
(see \RI, \S5, where an example of such functorial correspondence
has been constructed).

We now return to the geometric correspondence between $\JAB_{\GA}$ and $H$ (its detailed construction can be found in \RI, \S5.1).
This correspondence sends points of $\JAB_{\GA}$ to cycles of Calabi-Yau varieties.
More precisely, let $\ZA$ be a point in $\JAB_{\GA}$ and let $(z,[\alpha])$ be a point of the configuration $(Z,[\alpha])$. Then we attach to $(z,[\alpha])$ a particular
hyperplane section $H_{(z,[\alpha])}$ of $H$, which, according to Corollary \ref{h-sec}, is a Calabi-Yau variety. Our correspondence now is defined as the map
\begin{equation}\label{CY}
CY: \JAB_{\GA} \longrightarrow \PP(H^0(\OO_H (d)))
\end{equation}
which sends a point $\ZA$ to the divisor $CY \ZA = \displaystyle{\sum_{z\in Z} H_{(z,[\alpha])}}$ of degree $d =deg (Z)$ on $H$. So the essential part of constructing 
the map\footnote{The notation
$CY$ obviously stands for Calabi-Yau.} $CY$ is a definition of
the hyperplane $H_{(z,[\alpha])}$ for a point $(z,[\alpha]) \in (Z,[\alpha])$. 

For this we exploit the fact that the space
$H^0 (\OO_H (1)) = V^{\ast}$ is equipped with a distinguished basis $T,X_p,Y_p, (p=0,\ldots,\LG-2)$ defined in \S5.3. Thus one needs to produce constants
$T(z,[\alpha])$, $X_p (z,[\alpha])$, $Y_p (z,[\alpha])$, for $p=0,\ldots,\LG-2$, varying holomorphically with $(z,[\alpha])$, and then set $H_{(z,[\alpha])}$ to be the divisor corresponding
to the global section
$$
s(z,[\alpha]) = T(z,[\alpha]) T + \sum^{\LG-2}_{p=0} X_p (z,[\alpha]) X_p + \sum^{\LG-2}_{p=0} Y_p (z,[\alpha]) Y_p
$$
 of $\OO_H (1)$.

To calculate the values $T(z,[\alpha])$, $X_p (z,[\alpha])$, $Y_p (z,[\alpha])$, for $p=0,\ldots,\LG-2$, we use the path operator procedure described in \S5.2.
Namely, we take the delta-function $\delta_z$ on $Z$ supported at $z$, decompose it as in (\ref{d-f}) and ``propagate" the component
$\delta_{z,0} \in \HH^0 \ZA$ throughout the direct sum decomposition in (\ref{ordLG}) along the shortest (``geodesic") paths of the graph $G(\JABG)$
 joining the vertical level $0$ with the successive vertical levels  $p=1,\ldots,\LG -1$, and moving from left to right. In other words, we apply to
$\delta_{z,0}$ the operators $(D^{+} (\delta_{z,0}))^p$ to obtain the (right) moving string of functions
$$
\delta^{(p)}_{z,0} = (D^{+} (\delta_{z,0}))^p (\delta_{z,0})\,,
$$
for $p=0,\ldots, \LG-1$. Once we arrive to $\delta^{(\LG-1)}_{z,0} \in \HH^{\LG-1} \ZA$, we return this function back to $\HH^0 \ZA$ using the operators
$(D^{-} (\delta_{z,0}))^p$, for $p=1,\ldots,\LG-1$. This way we create the (left) moving string of functions
$$
\delta^{(\LG-1),(\LG-1-m)}_{z,0} = (D^{-} (\delta_{z,0}))^m (\delta^{(\LG-1)}_{z,0})\,.
$$
The desired constants are obtained essentially by evaluating
 all these functions at $z$. More precisely, we define
 \begin{eqnarray*}
 X_p(z,[\alpha]) =exp(\delta^{(p)}_z(z)),\,
 Y_p(z,[\alpha]) =exp(\delta^{(l-1),(p+1)}_z (z)),\,&\\
 T(z,[\alpha]) =exp(\delta^{(l-1),(0)}_z (z)),& 
 \end{eqnarray*}
 where $p=0,\ldots,\LG-2$.

The map $CY$ in (\ref{CY}) could be viewed as a variant of the classical Albanese map.
From this perspective the path operator technique together with the evaluation procedure described above
could be considered as  an analogue of integration of forms on the classical Jacobian along a path: the ``forms" in the case of
$\JA$ are sections of the sheaf $\HH^0 =\HT$ and the paths are the paths of the graph
$G(\JAB_{\GA})$.

The Calabi-Yau varieties attached to points $\ZA$ of $\JAB_{\GA}$ could be viewed as a geometric realization of the ``periods"
$\{\HH^p \ZA \}_{p=0,\ldots,\LG-1}$ attached to $\ZA \in \JAB_{\GA}$. One can say that our Jacobian allows to see how points
on a projective surface can ``open up" to become Calabi-Yau varieties. This is one of the striking predictions of quantum gravity according to the string theory and
it might be an indication that $\JA$ could be an effective mechanism to see the string theory behind the Hilbert schemes of $0$-dimensional subschemes.
\section{Example: Complete intersections}
We illustrate our general construction by considering complete intersections of
sufficiently ample divisors on a surface.
\\
\indent
Let $X$ be a smooth projective surface with irregularity
$q(X) = h^1(\OO_X) =0$. Fix a very ample line bundle $\OO_X (L)$
on $X$ and consider a $0$-dimensional complete intersection subscheme $Z$ of
two smooth irreducible members $C_1,C_2$ of the linear system $\mid L \mid$.
From the Koszul sequence 
$$
\xymatrix{
0\ar[r]& \OO_X (-2L)\ar[r]&\OO_X (-L) \oplus \OO_X (-L)\ar[r]& {\cal I}_Z 
\ar[r]&0 }
$$
for the ideal sheaf ${\cal I}_Z$ of $Z$ in $X$ we obtain
\begin{equation}\label{EZ=ls}
\EZ= H^0(\OO_X (L)) / {\bf C}\{s_1,s_2 \}\,,
\end{equation}
where $s_1$ and $s_2$ are sections of $H^0(\OO_X (L))$ corresponding to
$C_1$ and $C_2$ respectively. Thus 
$[Z] \in \stackrel{\circ}{\Gamma^r_d}(L)$, where $d=L^2$ and
$r =h^0(L) -3$. We always assume that
  $r \geq 1$ and set $\PP_Z$ to be the codimension
$2$ subspace of
 $\PP(H^0(\OO_X (L))^{\ast})$  spanned by $Z$ (we implicitly identify $Z$ with its image under the embedding given by
$\OO_X (L)$). 

By Lemma \ref{Htild}, 3), for a general choice of $\alpha \in \EZ$, there is a natural identification
of $\EZ$ with the subspace $\HT \ZA$ of $\HO Z)$. This identification can be described explicitly as follows.

Using
(\ref{EZ=ls}), a general choice of $\alpha \in \EZ$ amounts to choosing a section $s \in H^0(\OO_X (L))$
that does not vanish at any of the points of $Z$. Next we consider rational functions on $X$
given by quotients $\displaystyle{\frac{s^{\prime}}{s}}$, for $s^{\prime} \in H^0(\OO_X (L))$.
The restrictions of these functions to $Z$ are elements of $\HO Z)$ and they form the subspace
$\HT \ZA$ as $s^{\prime}$ varies in $ H^0(\OO_X (L))$. In other words we have
\begin{equation}\label{EZ=HT}
\HT \ZA = \{ \frac{s^{\prime}}{s}\,\, mod (I_Z) \mid s^{\prime} \in H^0(\OO_X (L)) \}\,,
\end{equation}
where $I_Z$ is the ideal of rational functions on $X$ vanishing on $Z$.
 This description implies that the morphism $\kappa_{\ZA}$ is simply the embedding
$$
Z \hookrightarrow \PP_Z
$$
defined by $\OO_X (L)$. In particular, the Hilbert function $P$ of
$Z$ in $\PP_Z$ determines the Hilbert function of 
the
filtration (\ref{H-filt}) over the points of $\XD$ corresponding to the
complete intersections of divisors in the linear system
$\left| L\right|$. 
Let $\JJ_Z$ be the ideal sheaf of $Z$ in $\PP_Z$. Then we have 
\begin{equation}\label{f1}
 P(k) = dim (\HT_{-k} ([Z], [\alpha])) = deg (Z) -h^1(\JJ_Z (k)) =
 L^2 -h^1(\JJ_Z (k))\,,
\end{equation}
where $\HT_{-k} ([Z], [\alpha])$ is the fibre of $\HT_{-k}$ in 
(\ref{H-filt}) at the point $([Z], [\alpha])$. In order to compute
$h^1(\JJ_Z (k))$ and, hence, $P(k)$ we make several assumptions on $L$:
\begin{enumerate}
\item[1)]
the line bundle $\OO_X (L)$ gives a projectively normal embedding of $X$,
\item[2)]
$h^2(\OO_X (kL)) = 0$, for all $k\geq 1$.
\end{enumerate}
(Observe that these assumptions are satisfied for a sufficiently high
multiple of any ample divisor on $X$). 
 
 These assumptions and straightforward computations (the details can be found in \RI, \S5.2) imply
\begin{equation}\label{P(k)}
P(k) = deg (Z) - h^2(\OO_X (kL)) + 2h^2(\OO_X ((k-1)L)) - 
h^2(\OO_X ((k-2)L))\,,
\end{equation}
for all $k\geq 0$. From 2) of the assumptions above it follows that
the filtration (\ref{H-filt}) at $([Z], [\alpha])$ 
 has the following form
\begin{equation}\label{f3}
0 = \HH_0 ([Z], [\alpha]) \subset \HH_{-1} ([Z], [\alpha]) \subset
\HH_{-2} ([Z], [\alpha])
\subset \HH_{-3} ([Z], [\alpha]) = \HO Z)\,.
\end{equation}

For $Z$ reduced, i.e. it is $L^2$ distinct points, and $[\alpha]$ in an appropriate non-empty Zariski open subset of
$\PP (\EZ)$,
 the filtration (\ref{f3}) splits to yield the following direct sum
decomposition.
\begin{equation}\label{f4}
\HO Z) = {\HH^0} ([Z], [\alpha]) \oplus {\HH^1} ([Z], [\alpha])
\oplus{\HH^2} ([Z], [\alpha])
\end{equation}
and we have 
\begin{alignat}{1} \label{f5}
dim({\HH^0} ([Z], [\alpha]))& =P(1) = 
h^0(\OO_X (L))-2\,,
\\
dim({\HH^1} ([Z], [\alpha]))& =
P(2) -P(1) = 
 L^2 - h^2(\OO_X) - h^0(\OO_X (L))+2\,,\nonumber  
\\
dim({\HH^2} ([Z], [\alpha]))& = 
P(3) -P(2) =
h^2(\OO_X)=h^0(\OO_X (K_X))\,.  \nonumber 
\end{alignat}

Let $\JAB$ be the part of $\JAA (X;L,L^2)$, where the summands of
the decomposition (\ref{f4}) have dimensions given by (\ref{f5}).
If the geometric genus $p_g(X) = h^0(\OO_X (K_X)) \neq0$, then
 Proposition \ref{H1} implies that the $(1,0)$-Dolbeault variety $H$ of $\JAB$
is a complete intersection of two quadrics 
$$
X_0 Y_0  = T^2,\,\,\, X_1 Y_1 = T^2
$$
in $\PP^4$ and its Lagrangian cycle $H_0$ is
the union of $4$ lines
$$
H_0 = \Pi_{\emptyset} \cup  \Pi_0 \cup  \Pi_1 \cup  \Pi\,,
$$
where $\Pi_{\emptyset} = \{X_0 =X_1 = 0\}$,
$\Pi_0= \{ X_1 = Y_0 =0\}$, $\Pi_1= \{ X_0 = Y_1 =0\}$ and
$\Pi =\{Y_0 = Y_1 = 0 \}$.

From Lemma \ref{singH} we obtain that $H$ is singular at $4$
points $(1:0:0:0),\,\,(0:1:0:0)\,\,,(0:0:1:0)\,\,,(1:0:0:0)$.
It is easy to see, either from projective or toric description of $H$,
 that resolving the singularities of $H$ we obtain
 $\PP^1 \times \PP^1$ blown-up at $4$ distinct points corresponding to
 the points of intersection of two reduced and reducible divisors
 $F_1 +F^{\prime}_1, \,\, F_2 +F^{\prime}_2$, where
 $F$ and $F^{\prime}$ are the divisor classes of the distinct rulings
 of
 $\PP^1 \times \PP^1$.
 \\
 \indent
 The Calabi-Yau cycle map $CY$ in Proposition \ref{CY} in this case
 sends points $([Z],[\alpha]) \in \JAB$ to a cycle of elliptic curves.
 More precisely, the construction in the proof of 
 Proposition \ref{CY} associates a smooth elliptic curve with every
 point $(z,[\alpha]) \in (Z,[\alpha])$.
 \\
 \indent
 Thus our construction implies that behind points on a smooth
 complex projective surface $X$ with $p_g \neq 0$ are ``hidden"
 elliptic curves. To ``reveal" them one has to make a point on $X$
 to be a part of a (reduced) complete intersection of curves
 in the linear system of a
  sufficiently high multiple of any ample divisor $L$ on $X$
 (observe that if $X$ is a K3-surface then taking {\it any} very ample
 $L$ will be enough).
 This could be viewed as an affirmative answer to the question
 of Nakajima about a possibility that
 ``elliptic curves are hidden in the Hilbert schemes" (see \cite{[N]}, p.2).
\section{Summary of {\bf Part I}}
The nonabelian Jacobian $\JA$ of a smooth complex projective surface $X$ is inspired by the classical theory of Jacobian of curves.
It is built as a natural scheme interpolating between the Hilbert scheme $\XD$ of subschemes of length
$d$ of $X$ and the stack
${\bf M}_X (2,L,d)$ of torsion free sheaves of rank $2$ on $X$ having the determinant $\OO_X (L)$ and the second Chern class (= number) $d$.

Viewing $\JA$ as a scheme over $\XD$ one uncovers various new structures of $\JA$. Namely, $\JA$ distinguishes certain admissible locally closed strata $\GA$
of $\XD$ which are labeled by weighted trivalent graphs $(G(\JABG),h^{\prime}_{\GA})$ as described in (\ref{G}). Each such weighted graph encapsulates

1) the period map intrinsically associated to $\JABG$ (see \S3.3),

2) the Fano toric variety parametrizing a natural family of Higgs structures attached to the sheaf $\FT \otimes \OO_{\JABG}$ (see \S4),

3) the sheaf of reductive Lie algebras {\BM $\LAGT$} on $\JABG$ (see \S3.5).
\\
\\
\noindent
The properties of the period map in 1), such as its injectivity (Torelli property), are related to algebro-geometric properties of the $0$-dimensional subschemes
of $X$ parametrized by $\GA$. The Fano toric varieties in 2) together with their hyperplane sections, which are (singular) Calabi-Yau varieties, open tantalizing possibilities
for new invariants for surfaces and vector bundles on them, based on the invariants of Fano toric varieties and Calabi-Yau varieties. This aspect of our Jacobian
also indicates at the intriguing ties with quantum gravity. A natural appearance of the sheaf of reductive Lie algebras in 3) opens a way to a systematic use of powerful methods
of the representation theory of Lie algebras/groups in the theory of projective surfaces. This aspect of $\JA$ is the subject of the second part of this survey.             
\newpage

\part{ Nonabelian Jacobian $\JA$ and representation theory}      

\section{Introduction}
 This part is entirely devoted to the sheaf of Lie algebras {\boldmath$\GT_{\mbox{\unboldmath${\GA}$}}$} defined in \S3.5.
Our intention is to make the representation theory of {\boldmath$\GT_{\mbox{\unboldmath${\GA}$}}$} to work for the geometry of $X$.
In other words, we ask: what kind of geometric properties of $X$ can be seen in the Lie theoretic properties of {\boldmath$\GT_{\mbox{\unboldmath${\GA}$}}$}?   
Pursuing this line of inquiry our considerations are naturally divided into two parts:
\\
\\
1) establish a dictionary between the properties of the sheaf of (reductive) Lie algebras {\boldmath$\GT_{\mbox{\unboldmath${\GA}$}}$} attached to every admissible component\footnote
{see Definition \ref{csa}.}
$\GA$ in $C^r (L,d)$ and geometric properties of subschemes parametrized by $\GA$; 
\\
\\
2) use the representation theory of {\boldmath$\GT_{\mbox{\unboldmath${\GA}$}}$} to define interesting objects (e.g. sheaves, complexes of sheaves) which can serve as new invariants of vector
bundles on $X$ as well as invariants of the surface itself.
\\
\\
\indent
For the first part we are able to uncover 
\\
\\
a) a precise relationship between the center of the reductive Lie algebras {\boldmath$\GT_{\mbox{\unboldmath${\GA}$}}$} and canonical decompositions
of configurations parametrized by $\GA$,
\\
b) how to use particular ${\bf \goth sl}_2$-subalgebras of our reductive Lie algebras to gain insight into geometry of configurations of points of $X$.
\\
\\
\indent
For the part 2) we show that the sheaf of reductive Lie algebras {\boldmath$\GT_{\mbox{\unboldmath${\GA}$}}$} attaches to $X$  a distinguished collection of
\\
\\
c) representations of symmetric groups and
\\
d) perverse sheaves of representation theoretic origin on the Hilbert scheme $\XD$.
\\
\\
 \indent
The appearance of the categories of the representations of symmetric groups in c) and perverse sheaves in d) 
comes from the fact that our Jacobian connects in a natural way to such fundamental objects in geometric representation theory as
Springer resolution of the nilpotent cone of simple Lie algebras (of type ${\bf A_n}$), Springer fibres, loop algebras and Infinite Grassmannians. 
Before proceeding with a more detailed account of these results, we would like to discuss how our approach fits into the landscape of some recent advances in the theory of Hilbert schemes
of points on surfaces.

The last 15 years saw some profound relations emerging between the representation theory and the Hilbert schemes of points on algebraic surfaces. One of them is undoubtedly comes
from the works of Grojnowski and Nakajima which shows that the direct sum of the cohomology rings
$\displaystyle{\bigoplus^{\infty}_{d=0} H^{\ast} (\XD,{\bf Q})}$ of {\it all} Hilbert schemes of a projective surface $X$ is an irreducible representation of certain affine algebras
(see \cite{[N]} for details and references). To get a feeling for the depth of this result it is enough to mention that as a corollary one obtains 
the formula of G\"ottsche, \cite{[Go]}, for the generating series of the Poincar\'e polynomials of the Hilbert schemes of points of $X$ as the character of the above mentioned representation. 
What is remarkable about G\"ottsche's formula is that each Poincar\'e polynomial has no nice expression but, when put together into the generating series, they form a beautiful expression
with a lot of structure suggesting its relation to the representation theory of affine Lie algebras. The result of Grojnowski and Nakajima establishes such a relation and
reflects the same phenomenon on the level of the cohomology rings - by itself each ring
$H^{\ast} (\XD,{\bf Q})$ has no representation theoretic pattern but when put together such a pattern appears. This is achieved by considering certain natural correspondences
between the Hilbert schemes $\XD$ and $X^{[d+i]}$, for every $d\geq0$ and every $i \in {\bf Z} \setminus \{0\}$. These correspondences put together give rise to operators of degree $i$
of
$\displaystyle{\bigoplus^{\infty}_{d=0} H^{\ast} (\XD,{\bf Q})}$. These are Nakajima operators $c[i]$, which are indexed by elements $c$ of the cohomology ring
$H^{\ast} (X,{\bf Q})$ of $X$ and integers $i$, the degree of those operators with respect to the grading of 
$\displaystyle{\bigoplus^{\infty}_{d=0} H^{\ast} (\XD,{\bf Q})}$.  Though it remains a mystery why these operators have the commutator relations of the Heisenberg/Clifford Lie
algebras, on the the intuitive level one can say, that by putting all Hilbert schemes together, one achieves a ``horizontal"\footnote{the adjective `horizontal' is used because we think of
all Hilbert schemes being put on the same horizontal line.}   dynamics between Hilbert schemes of different degrees.
And this ``horizontal" dynamics yields representation theoretic patterns on the level of the cohmology ring
$\displaystyle{\bigoplus^{\infty}_{d=0} H^{\ast} (\XD,{\bf Q})}$. 
All this suggests some kind of quantum phenomena and that, perhaps, the full understanding might come from new mathematics related to the string theory.

Heuristically,  the representation theoretic pattern in our story, i.e. the sheaf of Lie algebras {\boldmath$\GT_{\mbox{\unboldmath${\GA}$}}$}, also comes as a result
of certain dynamics of Hilbert schemes. But in our case the dynamics is ``vertical"\footnote{see the footnote on page \pageref{v} for the meaning of `vertical'.} and the whole construction is permeated with the string theoretic analogies. 
Namely, the Jacobian $\JA$ can be thought of as a thickening of the Hilbert scheme
$\XD$, where a point $[Z] \in \XD$ is considered together with its ``internal moduli" parameter $[\alpha] \in \PP(\EZ)$.\footnote{recall that $\PP(\EZ)$ is the fibre over $[Z]$ of the projection
$\pi:\JA \longrightarrow \XD$ in (\ref{pi-h}).} It is this parameter that turns $Z$ into a dynamical object.\footnote{for this of course one needs
the parameter $[\alpha] $ to have sufficient degree of freedom, which is the assumption $dim(\PP(\EZ))\geq 1$.} This dynamics is observed
not on $Z$ itself but on the cohomology ring $\HO Z)$ of $Z$. As we explained in \S3.2, for ``good"\footnote{i.e. $\ZA \in \JABG$, where for some admissible component 
$\GA$ in $C^r (L,d)$; the notation $\JABG$ is found in (\ref{pi-JAB}).} points $\ZA$ of $\JA$ one has a distinguished direct sum decomposition
\begin{equation}\label{ord-ZA-II}
\HO Z) = \bigoplus^{\LG}_{p=0} \HH^p \ZA
\end{equation}
which is the fibre at $\ZA$ of the sheaf decomposition in (\ref{ordFT}). Thus by going over to $\PP (\EZ)$, the fibre of the Jacobian $\JA$ over $[Z]$,  one discovers that 
the cohomology ring $\HO Z)$,
until then 
 with no interesting structure, acquires a Hodge-like decomposition.
 Furthermore, the variation of this decomposition as $[\alpha]$ moves in $\JAB_Z$, a certain Zariski open subset
of $\PP (\EZ)$, is one of the manifestations of the ``vertical" dynamics of $Z$ that we have alluded to above.

It should be pointed out that the $0$-th summand $\HH^0 \ZA$ in (\ref{ord-ZA-II}) is of special importance. There are at least three reasons for that.
The first one is that $\HH^0 \ZA$ is the same as $\HT \ZA$ (see  (\ref{ordH-1})) and hence, by Lemma \ref{Htild}, can be identified with the group of extensions $\EZ$. 
This is noteworthy because it characterizes the elements of $\HH^0 \ZA$ as those functions on $Z$ (the data of local nature) which have a global geometric meaning, i.e.
the global extensions of the form (\ref{ext-seq}) and, in particular, the rank $2$ torsion free sheaves on $X$ sitting in the middle of those exact sequences.

The second reason is that $\HH^0 \ZA =\HT \ZA$ can be viewed as a linear system on $Z$ and hence connects to the geometric properties of $Z$. In particular, it defines the morphism
$\kappa_{\ZA}$ which we encountered in (\ref{kappaZ}).
 
The third reason is contained in the identification of the quotient
$ \HH^0 \ZA / \CC \{1_Z \} $
with the vertical tangent space of $\JA$ at $\ZA$ (see Proposition \ref{H0-vertT}).

Thus $\HH^0 \ZA$ connects to the global geometry of $X$, sees the ``linear" part of certain projective properties of $Z$ and relates to the `vertical' infinitesimal variation of 
of the decomposition in (\ref{ord-ZA-II}). This very space also labels the operators generating the Lie algebra
  ${\bf \tilde{\goth g}} \ZA$, the fibre of {\boldmath$\GT_{\mbox{\unboldmath${\GA}$}}$} at $\ZA$. More precisely,  ${\bf \tilde{\goth g}} \ZA$ is defined
as the Lie subalgebra of $End(\HO Z))$ generated by the operators 
$D^{+} (t)$, $D^{-} (t)$ and $D^{0} (t)$ in the triangular decomposition in (\ref{d-Dt}), as $t$ ranges through the space $\HH^0 \ZA$.
Observe that the cohomolgy $\HO Z)$, in view of (\ref{ord-ZA-II}), is now graded and the operators $D^{\pm} (t)$ and $D^{0} (t)$, the generators of 
${\bf \tilde{\goth g}} \ZA$, are graded as well: they are of degree $\pm 1$ and $0$ respectively.

It is interesting to observe that our Lie algebra ${\bf \tilde{\goth g}} \ZA$, as in the case of Nakajima, is generated by operators on a graded cohomology ring in (\ref{ord-ZA-II}) and these
operators are labeled by vectors of a particular summand of the grading as well as their operator degree. But of course the Lie algebras
${\bf \tilde{\goth g}} \ZA$ are purely local, i.e. linked to points on {\it each} Hilbert scheme $\XD$, for appropriate values of $d$. The sheaf of Lie algebras
{\boldmath$\GT_{\mbox{\unboldmath${\GA}$}}$} is their globalization and, again, is related to each Hilbert scheme $\XD$, for appropriate values of $d$.
We expect that the utility of these Lie algebras lies precisely in their capacity to reveal some geometric properties of the underlying $0$-dimensional subschemes.
The results of \cite{[R2]} summarized below justify this thinking.

\section{The center of the Lie algebra ${\bf \tilde{\goth g}} \ZA$ and geometry of $Z$}
We determine the Lie algebras ${\bf \tilde{\goth g}} \ZA$ attached to points of the stratum $\JABG$, for each admissible component in $\CS$.
 In particular, their decomposition into the sum of simple Lie algebras is related to the geometry of the configurations parametrized by $\GAB$. More precisely, we show
\begin{thm} \label{tcd}
Let $\ZA$ be a point of $\JABG$, for some admissible component $\GA$ in $\CS$. Then $Z$ decomposes into the disjoint union
 \begin{equation}\label{cd}
Z = \bigcup_{i=1}^{\nu} Z^{(i)}\,,
\end {equation}
where $\nu$ is the dimension of the center of the Lie algebra ${\bf \tilde{\goth g}} \ZA$.
Furthermore, the Lie algebra ${\bf \tilde{\goth g}} \ZA$ and hence, its center act on the subspace $\HT_{-\LG} \ZA$ of the filtration
of $\HO Z)$ in (\ref{filtHT-JG}). This action of the center determines the weight decomposition
$$
\HT_{-\LG} \ZA =\bigoplus^{\nu}_{i=1} V_i \ZA
$$
which possesses the following properties:
\begin{enumerate}
\item[1)]
$\HO{Z^{(i)}}) \cong V_i \ZA \cdot \HO{Z})$,
\item[2)]
 one has a natural isomorphism
\begin{equation}\label{Lid}
{\bf  \tilde{\goth g}} \ZA \cong \bigoplus^{\nu}_{i=1} {\bf \goth gl}(V_i \ZA)\,.
\end {equation}
\end{enumerate}
\end{thm}
This result establishes a precise dictionary between the decomposition of the Lie algebra
${\bf  \tilde{\goth g}} \ZA$ into the direct sum of matrix algebras and the geometric decomposition of $Z$ into the disjoint union of subschemes
in (\ref{cd}).

It should be recalled that in discussing the differential of the period map $p_{\GA}$ on page \pageref{p-Zdec}, we remarked that the decomposition similar to the one in (\ref{cd})
is a consequence of the failure of this differential to be injective at a point $\ZA$ of $\JABG$. In fact, this is not a coincidence since
 the Lie algebra ${\bf  \tilde{\goth g}} \ZA$ also controls the properties of the derivative of the period map $p_{\GA}$ at $\ZA \in \JABG$.
\begin{thm}\label{infTor}
The derivative of the period map $p_{\GA}$ defined in (\ref{p}) is injective precisely at the points $\ZA$ of $\JABG$ for which
${\bf  \tilde{\goth{g}}} \ZA \cong {\bf \goth gl}_{d^{\prime}} (\CC)$, where $d^{\prime} =dim(\HT_{-\LG} \ZA)$.
\end{thm}
This is a version of the Infinitesimal Torelli Theorem for $\JA$ phrased in terms of the Lie algebraic properties of the sheaf
{\boldmath$\GT_{\mbox{\unboldmath${\GA}$}}$}.  Thus in our story the Infinitesimal Torelli property
(the injectivity of the differential of the period map $p_{\GA}$) has a precise geometric meaning:
it fails exactly when the decomposition (\ref{cd}) is non-trivial.

Let 
{\BM
$\CG$ be the center of $\LAGT$ and let $\LAG =[\LAGT,\LAGT]$ be its semisimple part. Then we have the structure decomposition 
\BEN\label{d-st}
\LAGT =\CG \oplus \LAG\,.
\EEN}
  The admissible components $\GA$ over which the period map $p_{\GA}$ itself is injective, i.e. the Torelli property holds for $p_{\GA}$,
can also be characterized by the properties of the sheaf  {\boldmath$\GT_{\mbox{\unboldmath${\GA}$}}$} and its structure decomposition. 
\begin{thm}\label{T=rk}
Let $\GA$ be a component in $\CSA$. Then the following statements are equivalent
\begin{enumerate}
\item[1)]
the Infinitesimal Torelli property holds for $p_{\GA}$
\item[2)]
the Torelli property holds for $p_{\GA}$
\item[3)]
the rank of 
{\BM
$\CG$}
is one. More precisely, 
{\BM
$\LAGT$}$=\ENDO (\HT_{-\LG} \otimes \OO_{\JABG})$,
its center 
{\BM
$\CG$}$\cong \OO_{\JABG}$ and its semisimple part 
{\BM
$\LAG$}$={\bf \goth sl} (\HT_{-\LG} \otimes \OO_{\JABG})$,
the sheaf of germs of traceless endomorphisms of 
$\HT_{-\LG} \otimes \OO_{\JABG}$.
\end{enumerate}
The components of $\CSA$ subject to (one of) these equivalent conditions is called simple.\footnote{if a component $\GA \in \CSA$ is not simple, then
the covering map 
$p_2 : \ZD_{\GAB} \longrightarrow \GAB$ admits a canonical non-trivial factorization determined by the characters of the action of the center
{\BM$\CG$} on the sheaf $\HT_{-\LG} \otimes \OO_{\JABG}$ (for more details see \cite{[R2]}, Theorem 3.26, and the discussion leading to it).}
\end{thm}

These results constitute a semisimple aspect of the representation theory of ${\bf  \tilde{\goth g}} \ZA$ in the sense that it takes into account the action of the center of 
${\bf  \tilde{\goth g}} \ZA$ (composed of semisimple elements) on the space $\HO Z)$. 
There are also nilpotent aspects  which are much more involved.

\section{Nilpotent aspects of ${\bf  \tilde{\goth g}} \ZA$}
{\bf 10.1. The maps $D_{\ZA}^{\pm}$.} The nilpotent aspect comes from the very definition of ${\bf  \tilde{\goth g}} \ZA$. Recall, that for every $t \in \HT \ZA$
 the multiplication operator $D(t)$ admits the triangular decomposition
\BEN\label{d-Dt-II}
D(t) =D^{+} (t) + D^{0} (t) + D^{-} (t)
\EEN
discussed in \S3.5 (see (\ref{d-Dt})). The operators $D^{\pm}(t)$ have degree $\pm 1$ with respect to the grading in
(\ref{ord-ZA-II}). In particular, these operators are nilpotent. Furthermore, if we  
let ${\bf \goth g} \ZA$  to be the semisimple part of ${\bf  \tilde{\goth g}} \ZA$, then $D^{\pm}(t) \in {\bf \goth g} \ZA$.

The algebraic meaning of these operators is easy to understand: if $x \in \HO Z)$ is of pure degree $p$, i.e. $x$ is a vector in $\HH^p \ZA$,
then $D^{+} (t)$ (resp. $D^{-} (t)$) records the part of the product $t\cdot x$ which goes over to the summand
$\HH^{p+1} \ZA$ (resp. $\HH^{p-1} \ZA$)\footnote{this is closely related to the formula for the relative differential of the period map in (\ref{diff-formula}); the relation is spelled out in
 \cite{[R2]}, Lemma 4.7.} . In particular, it follows that 
$D^{\pm} (t) =0$, for $t$ in the subspace ${\CC}\{ 1_Z \}$ of constants on $Z$. This together with the identification
$$
\HT \ZA / {\CC} \{ 1_Z \} \cong T_{\pi} \ZA
$$
in Proposition \ref{H0-vertT} imply that
 we can attach a nilpotent element $D^{+}(v)$ (resp. $D^{-}(v)$) of ${\bf \goth g} \ZA$ to every vertical tangent vector $v$ 
 at a point $\ZA \in \JABG$. As $v$ runs through the space $T_{\pi} \ZA$ of the vertical tangent vectors of $\JABG$ at $\ZA$ we obtain the linear maps
\begin{equation}\label{D+}
D_{\ZA}^{\pm} : T_{\pi} \ZA \longrightarrow \mbox{\boldmath${\cal N}$} ({\bf \goth g} \ZA)
 \end{equation}
of the vertical tangent space $T_{\pi} \ZA$ into the nilpotent cone {\boldmath${\cal N}$}$({\bf \goth g} \ZA)$ of ${\bf \goth g} \ZA$.

From the well-known fact that
{\boldmath${\cal N}$}$ ({\bf \goth g} \ZA)$ is partitioned into a finite set of nilpotent orbits we deduce that
the map\footnote{we discuss here the components $D^{+}$ of the decomposition (\ref{d-Dt-II}) but, of course, analogous results hold for $D^{-}$ as well.} $D_{\ZA}^{+} $
 assigns to $\ZA$ a finite collection of nilpotent orbits of {\boldmath${\cal N}$}$({\bf \goth g} \ZA)$.
These are the orbits intersecting the image of $D_{\ZA}^{+} $. Varying $\ZA$ in $\JABG$, for {\it simple} components $\GA \in \CSA$, we deduce the following.
\begin{thm}\label{nilo}
The Jacobian $\JA$ determines a finite collection
${\cal V}$ of quasi-projective subvarieties\footnote{${\cal V}$ is the collection of all simple admissible components of $\CS$ as $r$ ranges through all possible values of $r$.}
 of $\XD$ such that every $\Gamma \in {\cal V}$ determines a finite collection $O(\Gamma )$ of nilpotent orbits in 
${\bf \goth sl}_{d^{\prime}_{\Gamma}} (\CC)$, where  $d_{\Gamma}^{\prime} =rk (\HT_{-\LG} \otimes \OO_{\JABG})$.
\end{thm}
 Recalling that nilpotent orbits in ${\bf \goth sl}_n (\CC)$ are parametrized by the set of partitions $P_n$ of $n$, the above  result can be rephrased
by saying that every $\Gamma$ in ${\cal V}$ distinguishes a finite collection $P(\Gamma)$ of partitions of $d_{\Gamma}^{\prime}$.
Since partitions of $n$ also parametrize isomorphism classes of irreducible representations of the symmetric group
$S_n$ we obtain the following equivalent version of Theorem \ref{nilo}.
\begin{thm}\label{irrSn}
The Jacobian $\JA$ determines a finite collection
${\cal V}$ of quasi-projective subvarieties of $\XD$ such that every $\Gamma \in {\cal V}$
 determines a finite collection $R_{d^{\prime}_{\Gamma}}(\Gamma )$ of irreducible representations of the
symmetric group $S_{d^{\prime}_{\Gamma}}$, where  $d^{\prime}_{\Gamma}$ is as in Theorem \ref{nilo}.
\end{thm}

One way to express this result is by saying that the Jacobian $\JA$ elevates a single topological invariant $d$, the degree of the second Chern
class of sheaves parametrized by certain subvarieties of $\JA$, to the level of modules of symmetric groups.
Thus our Jacobian gives rise to new invariants with values in the categories of modules of symmetric groups.

But there is more to it. The partitions distinguished by $\JA$ contain a great deal of geometry of subschemes parametrized by $ \Gamma$'s in 
Theorem \ref{nilo}.
\\
\\
{\bf 10.2. $\goth{sl}_2$-subalgebras and equations of $\ZA$.} In down to earth terms, one can say that the partitions picked out by points $\ZA$ of $\JABG$, for $\GA \in {\cal V}$, yield equations defining
the image of $Z$ under certain morphisms\footnote{those are two morphisms of $Z$ naturally attached to our constructions: one is the morphism
$\kappa_{\ZA}$ in (\ref{kappaZ}) and the other one is given by the adjoint linear system $\left| L+K_X \right|$.} into appropriate projective spaces. The process of obtaining those 
equations is somewhat evocative
of the classical method of Petri (see \cite{[M]} for an overview). However, the essential ingredient in our approach is representation theoretic. 

It turns on the use of ${\bf \goth sl}_2$-subalgebras of ${\bf \goth g} \ZA$ associated to the nilpotent elements
$D_{\ZA}^{\pm} (v)$, the values of the maps $D_{\ZA}^{\pm}$ in (\ref{D+}). The operators $D_{\ZA}^{\pm} (v)$ in our considerations play the role of the operator\footnote
{the operator defined by the exterior product with a K\"ahler form; the standard notation ``L" for this operator, hopefully will not be confused with our use of the same notation.}
    ``L" in the Lefschetz decomposition in the Hodge theory. Completing $D_{\ZA}^{+} (v)$ (resp. $D_{\ZA}^{-} (v)$) to an ${\bf \goth sl}_2$-subalgebra of 
${\bf \goth g} \ZA$ in an appropriate way and considering the representation of this ${\bf \goth sl}_2$-triple on 
$\HO Z)$ gives us a sort of Lefschetz decomposition of $\HO Z)$. This combined with the orthogonal decomposition in (\ref{ord-ZA-II}) 
yields a bigrading of $\HO Z)$ thus revealing a much finer structure than the initial grading (\ref{ord-ZA-II}).
Once this bigrading is in place writing down the equations defining $Z$ in a certain projective space is rather straightforward.\footnote{to be precise, one obtains the equations defining the
{\it image} of $Z$ under the morphism $\kappa_{\ZA}$ in (\ref{kappaZ}) or under the morphism defined by the adjoint linear system $\left| L+K_X \right|$.}

 {\bf Example 1.} In this example we examine our constructions over
 the first non-trivial stratum $ \stackrel{\circ}{\GA^1_d}$ of the stratification of $\XD$ in (\ref{strat-Gam}). We assume that $d\geq 4$ and that there is a configuration of $d$ points $Z$,
such that
$[Z] \in \stackrel{\circ}{\GA^1_d}$.
 In addition, it will be assumed that $Z$ is in general position with respect to the adjoint linear system $\left| L+K_X \right|$.
This implies that the adjoint linear system  $\left| L+K_X \right|$ defines an embedding
\begin{equation}\label{Z-adj}
 Z \hookrightarrow \PP \big{(} (H^0(\OO_X ( L+K_X))/ H^0 (\ID_Z ( L+K_X)))^{\ast} \big{)} = \PP^{d-3}
\end{equation}
and any $d-2$ points of $Z$ span the ambient projective space $ \PP^{d-3}$ (we tacitly identify $Z$ and its image in $ \PP^{d-3}$).

A well-known result from the classical algebraic geometry says that $Z$ lies on a rational normal curve in $ \PP^{d-3}$ (see e.g. \cite{[G-H]}).
Let us show how this fact can be derived from our representation theoretic approach.
In what follows we assume $d\geq 5$, since for $d=4$ the assertion is obvious.

We begin by computing the filtration $\HT_{\bullet} \ZA$. From Lemma \ref{Htild}, 3), it follows that
$\HT_{-1} \ZA = \HT \ZA$ is two-dimensional.  This and the assumption that $Z$ is in general position with respect to the adjoint linear system $\left| L+K_X \right|$ 
implies that the morphism  $ \kappa_{\ZA} : Z \longrightarrow \PP ( \HT \ZA^{\ast}) = \PP^1$
in (\ref{kappaZ})
is an embedding as well. From this it follows that the filtration 
$$
0  \subset \CC\{1_Z \} \subset \HT_{-1} \ZA \subset \HT_{-2} \ZA \subset \cdots \subset \HT_{d-1} \ZA =\HO Z)
$$
is a maximal ladder. Thus the orthogonal decomposition in (\ref{ord-ZA-II}) has the following form
$$
\HO Z) =\bigoplus^{d-2}_{p=0} \HH^p \ZA
$$
and the Hilbert vector of $Z$ is $h_Z =(2,\overbrace{1,\ldots,1}^{\scriptscriptstyle{(d-2)-times}}\!,0)$.
 
The subspace $\FI^1 \ZA = \bigoplus^{d-2}_{p=1} \HH^p \ZA$, orthogonal to $\HT_{-1} \ZA$, can be naturally identified with the
space $H^0(\OO_X ( L+K_X))/ H^0 (\ID_Z ( L+K_X)$.\footnote{this follows from the dual version of the filtration $\HT_{\bullet}$ discussed on page \pageref{alt}.}

 We are going to write down a particular basis of $\FI^1 \ZA$ which viewed as
a basis of $H^0(\OO_X ( L+K_X))/ H^0 (\ID_Z ( L+K_X)$, will produce quadratic equations of $Z$ in $\PP^{d-3}$.

For this choose 
 $x \in  \HT_{-1} \ZA$ which is not in the subspace $ \CC\{1_Z \}$ of constants and consider the operator $D (x)$ of the multiplication by $x$.
We are interested in its action on $\FI^1 \ZA $.
So start with a generator $v_{d-2}$ of $\HH^{d-2} \ZA$ and apply to it powers of  $D (x)$ to produce vectors
\begin{equation}\label{vect}
v_{d-2-k} = x^k v_{d-2},
\end{equation}
 for $k=0,1,\ldots,d-3$. The main observation is that the vectors 
$\{v_{d-2},\ldots,v_1 \} $ form a basis of $\FI^1 \ZA$. This is the point where the nilpotent operator $D^{-} (x)$ comes into play. Namely, one uses the congruence
\begin{equation}
v_{d-2-k} = x^k v_{d-2} \equiv ( D^{-} (x))^k ( v_{d-2} )\,\,\, mod \{ \FI^{d-1-k} \ZA \}
\end{equation}
and the fact that $( D^{-} (x))^k ( v_{d-2} )$ is a generator\footnote{this fact comes from the strong
 form of the injectivity of the differential
$d_{\pi} (p_{\GA})$ at $\ZA$ which says that $d_{\pi} (p_{\GA})$ in (\ref{Grif-trans}) is injective
 component wise, see \cite{[R2]}, Proposition 4.21.} of $\HH^{d-2-k} \ZA$, for every $k=0,\ldots,d-3$.

We now think of $\{v_{d-2},\ldots,v_1 \} $ as basis  of linear forms on $\PP^{d-3}$ and arrange them in the following $2\times (d-3)$ matrix 
\BEN\label{matr-Z}
Q =\left( \begin{array}{ccc}
          v_{d-2} &\cdots&v_2 \\
           v_{d-3} &\cdots&v_1
               \end{array} \right)
\EEN
The $2\times 2$ minors of this matrix give us $\binom{d-3}{2}$ quadrics in $\PP^{d-3}$ cutting out a rational
normal curve. Furthermore, on $Z$ we have the relations (\ref{vect}). This implies that all of the
 above quadrics pass through $Z$.  Hence $Z$ lies on a rational normal curve.

 But actually we have more than just recovering an old result, since the rational curve
itself acquires now a new meaning - it can be identified with (the closure of) the image of period map associated to the filtration
$$
\FI^1 \ZA \supset \FI^2 \ZA  \supset \cdots  \FI^{d-2} \ZA  \supset  \FI^{d-1} \ZA  =0, 
$$
where  $ \FI^i =\bigoplus^{d-2}_{p=i} \HH^p \ZA$ and $[\alpha]$ varies in a suitable open subset of $\PP (\EZ) = \PP^1$.

From the representation theoretic point of view this example can be expressed by saying that $D^{-} (x)$, 
viewed as a nilpotent operator on the space $\FI^1 \ZA $,
is a principle nilponent element of ${\bf \goth{sl}} (\FI^1 \ZA) $ and it corresponds to the partition
$(d-2)$ of $d-2 =dim(\FI^1 \ZA)$. This classical example now admits a straightforward generalization.
Namely, if the partition corresponding to $D^{-} (x)$, viewed as a nilpotent operator on the space $\FI^1 \ZA $,
has $m$ rows\footnote{we think of a partition as the Young diagram associated to it.}, then $Z$ lies on the rational normal scroll of dimension $m$ and the geometry of this scroll is completely
determined by the shape of that partition (see \cite{[R2]}, \S9.6, for details).

{\bf Example 2.} In this example we apply our considerations to configurations which are complete intersections on a $K3$-surfaces.
In particular we give a complete set of explicit quadratic equations defining such configurations - the quadrics in question are all 
of rank
$\leq 4$. This in turn leads to recovering quadrics through canonical curve and those quadrics are much simpler than the ones 
obtained by Petri's method (see \cite{[M]}). Finally, one obtains a complete set of quadratic equations defining a $K3$-surface itself.

Let $X$ be a $K3$-surface and let $\OO_X (L)$ be a very ample line bundle on $X$. Consider a configuration $Z$ on $X$ which is a 
complete intersection of two smooth curves $C_1$ and $C_2$ in the linear system $\left| L \right|$. Let $\gamma_i \,(i=1,2)$ be sections of 
$H^0 (\OO_X (L))$ defining the curves $C_i \,(i=1,2)$, i.e. $C_i =(\gamma_i =0)$, for $i=1,2$. This is of course a special case of 
our considerations in 
\S6 and  
 as in (\ref{EZ=ls}) the space of extensions $\EZ$ can be identified as follows
\BEN\label{ext-ci}
\EZ = H^0 (\OO_X (L)) / {\CC\{\gamma_1, \gamma_2\}}.
\EEN
 From (\ref{f4}) the orthogonal decomposition of $\HO Z)$ at a point $\ZA$, for a general choice of 
$\alpha \in \EZ$, has the following form
\BEN\label{ord-ci}
\HO Z) = \HH^0 \ZA \oplus  \HH^1 \ZA \oplus  \HH^2 \ZA,
\EEN
where $dim  \HH^0 \ZA =\frac{L^2 }{2}$, $dim  \HH^1 \ZA =\frac{L^2 }{2} -1$ and $dim\HH^2 \ZA =1$.

From the identification in (\ref{ext-ci}) it follows that a choice of $\alpha $ in $\EZ$ amounts to
choosing a section $\gamma \in H^0 (\OO_X (L))$. In particular, a section $\gamma$ corresponding to a general $\alpha$
does not vanish at any point of $Z$ and as in (\ref{EZ=HT}) we have a natural identification
\begin{equation}\label{EZ=H0}
\HH^0 \ZA = \HT \ZA =\left\{ \left.\frac{\gamma^{\prime}}{\gamma}\right|_Z \mid \gamma^{\prime} \in H^0 (\OO_X (L)) \right\}\,,
\end{equation}
where $\left.\frac{\gamma^{\prime}}{\gamma}\right|_Z$ stands for the restriction of the rational function $\frac{\gamma^{\prime}}{\gamma}$ to $Z$.

Set $g=\displaystyle{\frac{L^2}{2} +1}$. This is the genus of smooth curves in the linear system $\left| L\right|$. Since we are aiming at writing quadratic equations defining
$X$, it will be assumed that $g\geq 4$. From now on we choose the complete intersection $Z$ in general
position in $\PP ((\EZ)^{\ast}) =\PP^{g-2}_Z$, the span of $Z$ in $\PP (H^0 (\OO_X (L))^{\ast}) =\PP^g$. Fix $g-1$ points $z_1, z_2,\ldots, z_{g-1}$ in $Z$ and view them 
as linear functionals on 
$\HH^0 \ZA $. The assumption that $Z$ is in general position means that $\{ z_1, z_2,\ldots, z_{g-1} \}$ is a basis of
$\HH^0 \ZA^{\ast}$. Let $\{x_1,x_2,\ldots, x_{g-1}\}$ be the basis of $\HH^0 \ZA $ dual to 
$\{ z_1, z_2,\ldots, z_{g-1} \} $, i.e.
\begin{equation}\label{dual-b}
x_i (z_j) =\delta_{ij},\,\,\forall i,j.
\end{equation}
We now consider the operators $D^{+} (x_i)$, for $i=1,\ldots,g-1$. In particular, we are interested in their restriction to $\HH^0 \ZA$, i.e. the homomorphisms
$$
D^{+} (x_i) : \HH^0 \ZA \longrightarrow \HH^1 \ZA
$$
The following properties of these maps are essential for writing down the quadratic equations for $Z$ in $\PP^{g-2}_Z$.
\begin{lem}\label{D+i}
 For every $i$, the homomorphism $D^{+} (x_i)$ is surjective and its kernel is the subspace $\CC\{1_Z \}$ of constants on $Z$.
\end{lem}
This implies that for every $i\neq j$ and for every $k$ there is an element $h^{\prime}_{ijk} \in \HH^0 \ZA $ such that
$$
D^{+} (x_i) (x_j) =D^{+} (x_k) (h^{\prime}_{ijk})
$$
This is a key relation which hides the quadratic equations defining $Z$ in $\PP_Z$.

To see the actual equations we
replace the operators $D^{+} (x_i) $ (resp. $D^{+} (x_k) $) by multiplication by $x_i$ (resp. $x_k$) to obtain a non-homogeneous  quadratic polynomial
\begin{equation}\label{quad}
q_{ijk}= x_i x_j - x_k h^{\prime}_{ijk} - m_{ijk} \in Sym^2 (\HH^0 \ZA)
\end{equation}
vanishing on $Z$ and where $ m_{ijk} \in \HH^0 \ZA $.
 We can be more precise about this linear term\footnote{Indeed, the quadratic part $x_i x_j - x_k h^{\prime}_{ijk}$ in (\ref{quad}) 
vanishes on the subset $\{z_1,\ldots,z_{g-1} \} \setminus \{z_k \}$.
Hence, $ m_{ijk}$ must vanish on those points. The general position condition assures that there is unique, up to a non-zero scalar, linear form on
$ \HH^0 \ZA^{\ast}$ vanishing on $\{z_1,\ldots,z_{g-1} \} \setminus \{z_k \}$ and this form is $x_k$.} - it must be a scalar multiple of $x_k$.
 Thus the polynomial in (\ref{quad}) has the form
 \begin{equation}\label{quad1}
q_{ijk}= x_i x_j - x_k ( h^{\prime}_{ijk}+c_{ijk}) \in Sym^2 (\HH^0 \ZA)\,,
\end{equation}
for some constant $c_{ijk}$. Furthermore, using the identification of $\HH^0 \ZA$ as the space of fractions in (\ref{EZ=H0}), we can complete the pencil $\{\gamma_1, \gamma_2 \}$ to a basis
$\{\gamma_1, \gamma_2, \gamma_3, \ldots,\gamma_{g+1}  \}$ of $H^0 (\OO_X (L))$ so that
$x_i =\displaystyle{\left.\frac{\gamma_{i+2}}{\gamma} \right|_Z}$, for every $i\in \{1,\ldots,g-1\}$.
Substituting these fractions in (\ref{quad1}) and clearing off the denominators we obtain homogeneous quadratic polynomials
$$
Q_{ijk} =\gamma_{i+2} \gamma_{j+2} - \gamma_{k+2} h_{ijk}
$$
 in $Sym^2\big(H^0 (\OO_X (L))\big)$ vanishing on $Z$.
Fixing $k=1$, and letting $i<j$ vary through the set $\{2,\ldots, g-1 \}$, we obtain
$\binom{g-2}{2}$ linearly independent quadratic polynomials vanishing on $Z$.  More precisely,
taking the restriction
$$
\overline{Q}_{ij1} =\overline{\gamma}_{i+2} \overline{\gamma}_{j+2} - \overline{\gamma}_{3} \overline{h}_{ij1}
$$
 of $Q_{ij1}$ to the codimension $2$ subspace $\PP_Z$ defines homogeneous quadratic equations of $Z$ in $\PP_Z$.

Let us choose a smooth curve in $\left| L\right|$, say $C_1 = (\gamma_1 =0)$, passing through $Z$ and denote
by $\ID_X$, $\ID_{C_1}$, and $\ID_Z$ the ideal sheaves respectively of $X$ in $\PP^g =\PP(H^0 (\OO_X (L))^{\ast})$, of its
hyperplane section $C_1$ in $\PP_{C_1} = \PP(\big{(}H^0 (\OO_X (L)) / \CC\{\gamma_1 \}\big{)}^{\ast})$ and of the complete intersection
$Z$ in $\PP_Z$. From the fact
$$
H^0 (\ID_X (2)) \cong H^0 (\ID_{C_1} (2)) \cong H^0 (\ID_Z (2))
$$
it follows that there exists a unique lifting of the quadratic polynomials $\overline{Q}_{ij1} \in H^0 (\ID_Z (2))$
to the polynomials
\begin{equation}\label{quad2}
\tilde{Q}_{ij1} = \gamma_{i} \gamma_{j} - \gamma_{3} h_{ij1} + a_{ij1} \gamma_1 + a_{ij2} \gamma_2, \,\,\mbox{for}\,\,4\leq i<j\leq g+1,
\end{equation}
in $H^0 (\ID_X (2)) $, where $a_{ij1}, a_{ij2} \in H^0 (\OO_X (L))$.  Thus we obtain a collection of $\binom{g-2}{2}$ linearly independent sections in $H^0 (\ID_X (2)) $. On the other hand
$$
h^0 (\ID_X (2)) = \binom{g+2}{2} - h^0 (\OO_X (2L))= \binom{g+2}{2} - (4g-2) =\binom{g-2}{2}.
$$
Hence the collection of quadratic polynomials
$\{\tilde{Q}_{ij1} \mid 4\leq i<j\leq g+1 \}$ in (\ref{quad2}) is a basis of $H^0 (\ID_X (2))$ and their restriction to the hyperplane
$\PP_{C_1}$ give the quadratic polynomials
\begin{equation}\label{quad3}
\tilde{Q^{\prime}}_{ij1} = \omega_{i} \omega_{j} - \omega_{3} \omega_{ij1} + \omega_2 \omega^{\prime}_{ij1}  , \,\,\mbox{for}\,\,4\leq i<j\leq g+1,
\end{equation}
where $\omega_k =\left.\gamma_k \right|_{C_1}$, for $2\leq k \leq g+1$, $\omega_{ij1} =\left.\gamma_{ij1} \right|_{C_1}$ and
$\omega^{\prime}_{ij1} = \left.a_{ij1} \right|_{C_1}$. Those are equations of a complete set of linearly independent quadrics cutting out the
canonical curve $C_1 \hookrightarrow \PP_{C_1}$.
\\
\\
\indent
The two {\it a priori} different facts of classical algebraic geometry considered above, appear now as manifestations of the same phenomenon - 
the orthogonal decomposition attached to configurations of points and the representation theory associated to this decomposition.
Thus our Jacobian comes along with a unifying representation theoretic method for treating various questions of geometry of surfaces. 
 
Of course, in general, the equations obtained by this representation theoretic method can be complicated and not very illuminating. 
What is essential in our approach is that this possibly very complicated set of equations is encoded in an appropriate 
${\bf \goth sl}_2$-decomposition of $\HO Z)$. This in turn can be neatly ``packaged" in the properties of the partitions singled out by the points
$\ZA$ of $\JABG$ ``polarized" by operators $D_{\ZA}^{\pm} (v)$, with $v$ varying in $T_{\pi} \ZA$ as in (\ref{D+}).

To summarize one can say that the nilpotent aspect of the representation theory of {\boldmath $\LAG$}, for $\GA \in \CSA$,  provides new geometric insights as well as 
new invariants of the representation theoretic nature. 
But this turns out to be only a part of the story. In fact we can go further: from $\JA$ to the category of perverse sheaves 
on $\XD$.
\\
\\
{\bf 10.3. From $\JA$ to perverse sheaves on $\XD$.}
The perverse sheaves discussed here arise from the pullback of the Springer fibration over
the nilpotent cone {\boldmath${\cal N}$}$({\bf \goth g} \ZA)$ of ${\bf \goth g} \ZA$ to the vertical tangent space $T_{\pi} \ZA$ via the map
$D^{+}_{\ZA}$ in (\ref{D+}). In fact, to simplify the matters, the cohomology ring of the Springer fibres gives rise to local systems on the smooth part $\GA^{sm}$ of simple
admissible components $\GA$ of $\XD$. Taking the associated Deligne-Goresky-MacPherson complex on the closure $\overline{\GA}$ of $\GA$ in $\XD$ and then extending
it by zero to the whole $\XD$ give us a distinguished collection of perverse sheaves on $\XD$.    
\begin{thm}\label{ps}
The Jacobian $\JA$ determines a finite collection 
{\boldmath${\cal P}$}$(X;L,d)$ of perverse sheaves on $\XD$.
These perverse sheaves are parametrized by pairs
$ (\Gamma,\lambda)$, where $\Gamma$ is a subvariety in ${\cal V}$ as in Theorem \ref{nilo} and $\lambda$ is a partition in $P(\Gamma)$.
\end{thm}
This result subsumes two previous theorems since  the perverse sheaves 
${\cal C}(\Gamma,\lambda)$ in {\boldmath${\cal P}$}$(X;L,d)$ have the following properties:
\\
a) ${\cal C}(\Gamma,\lambda)$ is the Intersection Cohomology complex
$IC(\Gamma^{sm}, {\cal L}_{\lambda})$ associated to a certain local system ${\cal L}_{\lambda}$ on $\Gamma^{sm}$.
\\
b) The local system ${\cal L}_{\lambda}$ corresponds to a representation 
\begin{equation}\label{repgf}
\rho_{\Gamma,\lambda} : \pi_1 ( \Gamma^{sm},[Z]) \longrightarrow Aut(H^{\bullet}(B_{\lambda}, {\bf C}))
\end{equation}
of the fundamental group $\pi_1 ( \Gamma^{sm},[Z])$ of $\Gamma^{sm}$ based at a point $[Z]\in \Gamma^{sm}$
and where 
$H^{\bullet}(B_{\lambda}, {\bf C})$ is the cohomology ring (with coefficients in $\bf C$) of a Springer fibre
$B_{\lambda}$ over the nilpotent orbit {\boldmath$O_{\mbox{\unboldmath$\lambda$}}$} of ${\bf \goth sl}_{d_{\Gamma}^{\prime}} (\CC)$
corresponding to the partition $\lambda$.
\\
c) The representation $\rho_{\Gamma,\lambda}$ admits the following factorization
\begin{equation}\label{Sprep}
\rho_{\Gamma,\lambda} : \pi_1 ( \Gamma^{sm},[Z]) \stackrel{\rho^{\prime}}{\longrightarrow} S_{d_{\Gamma}^{\prime}} \stackrel{sp_{\lambda}}{\longrightarrow} Aut(H^{\bullet}(B_{\lambda}, {\bf C})),
\end{equation}
where 
$S_{d_{\Gamma}^{\prime}} \stackrel{sp_{\lambda}}{\longrightarrow} Aut(H^{\bullet}(B_{\lambda}, {\bf C}))$ 
is the Springer representation of the Weyl group
$W=S_{d_{\Gamma}^{\prime}}$ of ${\bf \goth sl}_{d_{\Gamma}^{\prime}} (\bf C)$ on the cohomology of a fibre $B_{\lambda}$
of the Springer resolution
\begin{equation}\label{Spres}
\sigma: \mbox{\boldmath$\tilde{\cal N}$} \longrightarrow \mbox{\boldmath${\cal N}$}({\bf \goth sl}_{d_{\Gamma}^{\prime}} ({\bf C}))
\end{equation}
of the nilpotent cone $\mbox{\boldmath${\cal N}$}({\bf \goth sl}_{d_{\Gamma}^{\prime}} ({\bf C}))$ of ${\bf \goth sl}_{d_{\Gamma}^{\prime}} (\bf C)$
and where a fibre $B_{\lambda}$ is taken over the nilpotent orbit {\boldmath$O_{\mbox{\unboldmath$\lambda$}}$}
in $\mbox{\boldmath${\cal N}$}({\bf \goth sl}_{d_{\Gamma}^{\prime}} ({\bf C}))$.
\\
\\
\indent
Using the fact that the category of perverse sheaves is semisimple the collection {\boldmath${\cal P}$}$(X;L,d)$ gives rise to a distinguished
collection, denoted $\mbox{\BM$C$}(X;L,d)$, of {\it irreducible} perverse sheaves on $\XD$. This in turn defines the abelian category
${\cal A} (X;L,d)$ whose objects are isomorphic to finite direct sums of complexes of the form ${\cal C}[n]$, where 
${\cal C} \in \mbox{\BM$C$}(X;L,d)$ and $n\in {\bf Z}$. 

The construction of perverse sheaves outlined above is somewhat reminiscent of the construction of local systems on the classical Jacobian.
Recall, that if $J(C)$ is the Jacobian of a smooth projective curve $C$, then isomorphism classes of irreducible local systems on $J(C)$ are
parametrized by the group of characters $Hom (H_1 (J(C)), \CC^{\times})$. We suggest that the collection of irreducible perverse sheaves
$\mbox{\BM$C$}(X;L,d)$ could be viewed as a non-abelian analogue of the group of characters of the classical Jacobian, while
the abelian category ${\cal A} (X;L,d)$ could be envisaged as an analogue of the group-ring of $Hom (H_1 (J(C)), \CC^{\times})$.

Though objects of  ${\cal A} (X;L,d)$ are complexes of sheaves on the Hilbert scheme $\XD$, they really descend from $\JA$ and one of the ways to
remember this is the following
\begin{thm}\label{exp-int}
Let $\stackrel{\circ}{\JAA} (X;L,d) = \JA \setminus \mbox{\BM$\Theta$} (X;L,d)$ be the complement of the theta-divisor
$\mbox{\BM$\Theta$} (X;L,d)$ in $\JA$ and let
${\cal T}^{\ast}_{\stackrel{\circ}{\JAA} (X;L,d)/ {\XD}}$ be the sheaf of relative differential forms of 
$\stackrel{\circ}{\JAA} (X;L,d)$ over $\XD$. 
Then there is a natural map
$$
exp(\int): H^0 ({\cal T}^{\ast}_{\stackrel{\circ}{\JAA} (X;L,d)/ {\XD}} ) \longrightarrow  {\cal A} (X;L,d)\,.
$$
\end{thm}
The map in the above theorem could be viewed as a reincarnation of the classical map
$$
H^0 ({\cal T}^{\ast}_{J(C)} ) \longrightarrow Hom (H_1 (J(C)), \CC^{\times})\,,
$$
where ${\cal T}^{\ast}_{J(C)}$ is the cotangent bundle of $J(C)$. This map sends a holomorphic $1$-form $\omega$ on $J(C)$ to the exponential of the linear functional
$$
\int (\omega) : H_1 (J(C)) \longrightarrow \CC
$$
given by integrating $\omega$ over $1$-cycles on $J(C)$ (the notation `$exp(\int)$' in Theorem \ref{exp-int} is an allusion to this classical
map).
\\
\\
\indent
Relations of the Hilbert schemes of points of surfaces to partitions is not new. Notably, Haiman's work on the Macdonald positivity conjecture,
\cite{[Hai]}, makes an essential use of such a relation. The same goes for an appearance of perverse sheaves on $\XD$: the work of 
G\"{o}ttsche and Soergel, \cite{[Go-S]}, uses the decomposition theorem of \cite{[BBD]} for the direct image of the Intersection
cohomology complex $IC(\XD)$ under the Hilbert-Chow morphism to compute  the cohomology of Hilbert schemes. In both of these works the partitions appear from the outset because the authors exploit the points of the Hilbert scheme corresponding to the zero-dimensional subschemes $Z$ of $X$, where the points in $Z$ are allowed to collide according to the pattern determined by partitions. In our constructions it is essential to work over the open
part $Conf_d (X)$ of $\XD$ parametrizing configurations of $d$ distinct points of $X$. So there are no partitions seen on the level of the Hilbert scheme.
The partitions become visible only on the Jacobian $\JA$ via the Lie algebraic invariants attached to it. This again can be attributed to the phenomenon, already discussed in the
introduction to {\bf Part II}, that 
our constructions turn a configuration of distinct points with no interesting structure on it into a dynamical object. The dynamics here has a precise meaning: it is
 given by certain
linear operators acting on the space of complex valued functions on a configuration. In particular, the operators 
$D^{\pm}_{\ZA} (v)$ obtained as values of the morphisms $D^{\pm}_{\ZA}$ in (\ref{D+}) give rise to the ``propagations" and ``collisions" in the direct sum decomposition (\ref{ord-ZA-II}). 
This is not an actual, physical, collision of points in a configuration but rather algebro-geometric constraints for
a configuration to lie on hypersurfaces in the appropriate projective spaces. The partitions attached to the nilpotent operators $D^{\pm}_{\ZA} (v)$
can be viewed as combinatorial (or representation theoretic) measure of this phenomenon, while the perverse sheaves in Theorem \ref{ps} are its categorical expression.
\section{From $\JA$ to Affine Lie algebras}\label{AffLie}
By now we can say that the sheaf of Lie algebras  
{\boldmath$\LAGT$} establishes a solid bridge between the geometry of $X$ and the representation theory of reductive Lie algebras.
In this section we show that the construction of {\boldmath$\LAGT$} as well as of its semisimple part {\boldmath$\LAG$} is flexible enough to allow a natural appearance
of the loop analogues
 of these Lie algebras.
One can say that our Jacobian is capable of `seeing' loops on the Hilbert schemes of points on a projective surface. At this stage we do not know if the affine Lie algebras
in our constructions could be linked to the representations of the affine Lie algebras on the cohomology ring of all Hilbert schemes discovered by Grojnowski and 
Nakajima.\footnote
{see the discussion in \S8.} But since the appearance of the affine Lie algebras in our constructions is natural, it is plausible to expect that such a relation might exist.

It is clear, that formally we can replace the semisimple Lie algebra ${\bf \goth g}\ZA$ attached to a point $\ZA \in \JA$ by its loop Lie algebra
${\bf \goth g}\ZA [z^{-1},z]$, where $z$ is a formal variable. However, there is a more natural and explicit reason for appearance of loop
Lie algebras in our story. To explain this we recall that the Lie algebra
${\bf \goth g}\ZA$ is the semisimple part of the reductive Lie algebra ${\bf \tilde{\goth g}}\ZA$ that is defined as follows.

For every $h$ in the summand $\HH^0 \ZA$ of the decomposition (\ref{ord-ZA-II}) we consider the operator $D(h)$ of the multiplication by $h$ in the ring
$\HO Z)$. Decomposing this operator according to the direct sum in (\ref{ord-ZA-II}) yields a triangular decomposition
\begin{equation}\label{tri-d}
D(h)= D^{-}(h) +D^0(h) + D^{+}(h)
\end{equation}
which was discussed in \S3.5. It is quite natural and immediate to turn (\ref{tri-d}) into a loop
\begin{equation}\label{l-tri-d}
D(h,z)= z^{-1} D^{-}(h) +D^0(h) + zD^{+}(h),
\end{equation} 
where $z$ is a formal parameter. This natural one-parameter deformation of the multiplication in $\HO Z)$ is behind the following
loop version
of the map (\ref{D+}):
\begin{equation}\label{lD+}
LD^{+}_{\ZA} :\, \stackrel{\circ}{T}_{\pi} \ZA \longrightarrow Gr({\bf \goth g}\ZA)\,,
\end{equation}
where $Gr({\bf \goth g}\ZA)$ is the loop or Infinite Grassmannian of the semisimple Lie algebra ${\bf \goth g}\ZA$
and $\stackrel{\circ}{T}_{\pi} \ZA$ is an appropriate Zariski open subset of the vertical tangent space ${T}_{\pi} \ZA$ of $\JA$ at 
$\ZA$.
This gives the following `loop' version of Theorem \ref{nilo}
\begin{thm}\label{lo}
The Jacobian $\JA$ gives rise to a finite collection $\cal V$ (the same as in Theorem \ref{nilo}) of subvarieties $\Gamma$ of
$Conf_d (X)$.
Every such $\Gamma$ determines a finite collection 
$LO(\Gamma)$ of orbits of the Infinite Grassmannian
$Gr({\bf SL}_{d^{\prime}_{\Gamma}} ({\bf C}))$
of ${\bf SL}_{d^{\prime}_{\Gamma}}({\bf C})$, 
where $d^{\prime}_{\Gamma}$ is the same as in Theorem \ref{nilo}.
\end{thm}
Taking the Intersection Cohomology complexes
$IC({O})$ of the orbits $O$ in $LO(\Gamma)$, for every $\Gamma$ in $\cal V$, we pass to the category of perverse sheaves on 
$Gr({\bf SL}_{d^{\prime}_{\Gamma}} ({\bf C}))$. Now a result of Ginzburg, \cite{[Gi]}, and Mirkovi\v{c} and Vilonen, \cite{[M-V]},
which establishes an equivalence between the category of perverse sheaves (subject to a certain equivariance condition) on the Infinite Grassmannian $Gr({\bf G})$ of a semisimple
Lie group ${\bf G}$ and the category of finite dimensional representations of the Langlands dual group $\bf {}^L G$ of $\bf G$, gives a Langlands dual
version of Theorem \ref{nilo}.
\begin{thm}\label{LD}
For every subvariety $\Gamma$ in $\cal V$ in Theorem \ref{lo} the Jacobian
$\JA$ determines a finite collection ${}^L R(\Gamma)$ of irreducible representations of the Langlands dual
group
${}^L {\bf SL}_{d^{\prime}_{\Gamma}}({\bf C}) = {\bf PGL}_{d^{\prime}_{\Gamma}}({\bf C})$.
\end{thm}

In retrospect a connection of our Jacobian with the Langlands duality could have been foreseen. After all the nature of
$\JA$ as the moduli space of pairs $(\SEE)$ resembles the moduli space of pairs of Drinfeld in \cite{[Dr]}.
The fundamental difference is that the groups
${\bf SL}_{d^{\prime}_{\Gamma}}({\bf C})$ and their Langlands duals in our story have nothing to do with the structure group
($\bf GL_2 (C)$) of bundles parametrized by $\JA$. These groups rather reflect the geometric underpinnings of our construction
related to the Hilbert scheme $\XD$. Noting this difference, we also point out one of the key features of $\JA$:
\BEN\label{ff-dict}
\mbox{\it it transforms vertical vector fields on $\JA$ to perverse sheaves on $\JA$. } 
\EEN
By vertical vector fields in (\ref{ff-dict}) we mean sections of the relative tangent sheaf
${\cal T}_{\pi} ={\cal T}_{\JA /\XD}$ of $\JA$. In fact, (\ref{ff-dict}) is somewhat more general, since it applies to vertical vector fields supported on the locally closed
subset of $\JA$ of the form $\JABG$, for $\GA$ simple admissible component in $\CS$. This feature is essentially the map in Theorem \ref{exp-int} and it can be 
viewed as a ``tangent" version of Grothendieck's ``functions-faisceaux dictionnaire"\footnote{the correspondence
in (\ref{ff-dict}) also goes in the direction opposite to the Grothendieck's which assigns the function to a perverse sheaf.}
 which plays an important role in a reformulation of the classical, number theoretic, Langlands correspondence into the geometric one (see \cite{[Fr]},
 for an excellent introduction to the subject of the geometric Langlands program).

\section{$\JA$ and Langlands correspondence for surfaces}\label{Langlands}
A discussion of the subject of Langlands correspondence will take us too far a field.
So we just sketch here our point of view and indicate how the Jacobian $\JA$ could be used in this context.

One of our interests is to study correspondences of a smooth projective surface $X$. By this we mean subvarieties of (co)dimension $2$ in the
Cartesian product $X\times X$. The main idea of our approach consists of associating to such a subvariety $Y$ invariants/objects of the representation theoretic nature.
Of course, the ``dream" would be to have enough of such invariants/objects to recover $Y$ in $X\times X$. This is a kind of Tanakian philosophy and it is also reminiscent of the 
Langlands program to study the Galois group of an algebraic closure of a number field. In our case the role of the number field is played by the surface $X$ and an irreducible
correspondence $Y\subset X\times X$ is viewed, via the projections\footnote{we tacitly assume that projections induce finite (ramified) coverings $Y\longrightarrow X$.}
 on each factor of $X\times X$, as a finite algebraic extension of $X$. Thus we suggest to view correspondences of $X$ as (a part of) the Galois side of a hypothetical Langlands correspondence
for surfaces. 

What should be the automorphic side of such a correspondence? The lesson of the Geometric Langlands program for curves, \cite{[Fr]}, teaches us that it should be a suitable
subcategory of perverse sheaves on the moduli stack of vector bundles on $X$. So, roughly speaking, the Geometric Langlands correspondence for surfaces
 seeks to associate to an irreducible, non-trivial correspondence $Y$
in $X\times X$ a collection or a subcategory of perverse sheaves on a suitable stack of vector bundles on $X$. 

How does all this relates to the constructions presented so far? Our first remark concerns the universal scheme $\ZD$ and the covering morphism
\BEN\label{uc1}
p_2 : \ZD \longrightarrow \XD
\EEN
encountered in (\ref{uc}). We suggest to view $\ZD$ as a scheme {\it encompassing all} non-trivial correspondences of $X$ of degree $d$.\footnote{ this is the universality property of the morphism 
$p_2 : \ZD \longrightarrow \XD$, i.e. any non-trivial correspondence of $X$ of degree $d$ is obtained
 as the pullback of the universal family (\ref{uc1}).
More precisely, an irreducible correspondence $Y \subset X\times X$
is said to be of degree $d$ if one of the projections $X\times X \longrightarrow X$ induces a finite and flat morphism $p:Y \longrightarrow X$ of degree $d$. The universality of the Hilbert scheme
$\XD$ implies that there is a unique morphism $c:X\longrightarrow \XD$ such that $p:Y \longrightarrow X$ is the pullback of the universal family (\ref{uc1}) via the morphism $c$.} From this perspective
 the function field
${\bf K}(\ZD)$ of $\ZD$, an extension of degree $d$ of the function field ${\bf K}(\XD)$ of $\XD$, and its Galois group $G({\bf K}(\ZD) /{\bf K}(\XD)) =S_d$ 
appear as objects relevant to the Galois side of the Langlands correspondence for surfaces.

The second remark concerns the Jacobian $\JA$. It is a scheme interpolating\footnote{see the diagram (\ref{pi-h}).} 
 between the Hilbert scheme $\XD$ and the the moduli stack ${\bf M}_X (X,L,d)$ of torsion free sheaves of rank $2$
on $X$ with the fixed determinant $\OO_X (L)$ and the second Chern number $d$. So suitable perverse sheaves on $\JA$ could be viewed as an ``approximation" of the
automorphic side of the Langlads correspondence for surfaces.

From the above remarks it follows that establishing correspondences between the function field ${\bf K}(\ZD)$ and the symmetric group $S_d$ on the one side, and the category of perverse sheaves on
$\JA$ on the other side, could be viewed as a step toward Geometric Langlands correspondence for surfaces.
     
This is exactly what the dictionary (\ref{ff-dict}) does for us - it allows to attach perverse sheaves on $\JA$ to

(i) rational functions on the universal scheme $\ZD$ in (\ref{uc1}) and to

(ii) symmetric polynomials in $d$ indeterminates.
\\
\noindent
We explain the main points of these constructions in the following two subsections.
\\
\\
\noindent
{\bf 12.1. From ${\bf K}(\ZD)$ to perverse sheaves on $\JA$.} 
Let $f$ be a rational function on $\ZD$. According to (\ref{ff-dict}) it is enough to transform $f$ into a vertical vector field on $\JA$. More precisely, what we will show is how
to transform $f$ into a vertical vector field on the strata $\JABG$ (see (\ref{pi-JAB}) for notation) of $\JA$, for suitable simple admissible components\footnote{see Theorem \ref{T=rk} for definition.} $\GA \in \CS$.

Our starting point is to view $f$ as a rational section, call it $s_f$,
of the direct image $p_{2\ast} \OO_{\ZD}$, where $p_2$ is a morphism in (\ref{uc1}). Let us assume that there is a simple admissible component $\GA \in \CS$
intersecting the domain of the definition of $s_f$ non-trivially\footnote{this assumption holds for a general choice of $f$.}, i.e. the restriction $\left.s_f \right|_{\GA}$ is a rational section of 
$(p_{2\ast} \OO_{\ZD}) \otimes \OO_{\GA}$ which is regular on a non-empty Zariski open subset $\GA^{\prime}$ of $\GA$. By restricting further, we may assume 
that $s_f$ is regular on a non-empty Zariski open subset $\GA^{\prime\prime}$ contained in $\GAB$ (see (\ref{GAB}) for notation). Denote by $s^{\prime\prime}_f$
the restriction of $s_f$ to $\GA^{\prime\prime}$ and take its pullback
$\tilde{s}^{\prime\prime}_f =\pi^{\ast}s^{\prime\prime}_f$ to obtain a section of
$\FT \otimes \OO_{\pi^{-1} (\GA^{\prime\prime})} = \PIPO\otimes \OO_{\pi^{-1} (\GA^{\prime\prime})} $. Using the orthogonal decomposition of 
$\FT\otimes \OO_{\JABG}$ in (\ref{ordFT}) we write
 \BEN\label{sf-ord}
\tilde{s}^{\prime\prime}_f = \sum^{\LG}_{p=0} \tilde{s}^{\prime\prime}_{f,p},
\EEN
where $\tilde{s}^{\prime\prime}_{f,p} \in H^0 (\pi^{-1} (\GA^{\prime\prime}), \HH^p \otimes \OO_{\pi^{-1} (\GA^{\prime\prime})})$
is the $p$-th component of $ \tilde{s}^{\prime\prime}_f$ in the decomposition (\ref{ordFT}).  It is the component $\tilde{s}^{\prime\prime}_{f,0}$
that we are after. This is a section of 
$\HH^0 \otimes \OO_{\pi^{-1} (\GA^{\prime\prime})}$. We now recall the splitting
$$
\HH^0 \otimes \OO_{\pi^{-1} (\GA^{\prime\prime})} =\OO_{\pi^{-1} (\GA^{\prime\prime})} \oplus \HH^{\prime} \otimes \OO_{\pi^{-1} (\GA^{\prime\prime})}   
$$
from (\ref{f:splitHT}) and the identification of $\HH^{\prime} $ with the relative tangent bundle ${\cal T}_{\pi}$ of $\JABG$ in (\ref{Mpr}). 
These two facts imply that the projection of $\tilde{s}^{\prime\prime}_{f,0}$ onto $\HH^{\prime} \otimes \OO_{\pi^{-1} (\GA^{\prime\prime})}$
gives rise to a section of  ${\cal T}_{\pi}$ defined over the Zariski open subset $\pi^{-1} (\GA^{\prime\prime})$ of $\JABG$. Now we are in the position to
use (\ref{ff-dict}) to obtain a perverse sheaf on $\JA$.
\\
\\
{\bf 12.2. From symmetric polynomials to perverse sheaves on $\JA$.} Let $f=f(T_1,\ldots, T_d)$ be a symmetric polynomial in $d$ indeterminates $T_1,\ldots, T_d$.
We explain how to attach to $f$ a section $s^{\GA}_f$ of 
$\FT \otimes \OO_{\JABG}$, for every admissible component $\GA$ in $\CS$. Once such a section is constructed we repeat the argument in \S12.1, to
obtain the corresponding vertical vector field on $\JABG$.

To go from $f$ to a section $s^{\GA}_f$ of $\FT \otimes \OO_{\JABG}$ it will be enough to describe the values $s^{\GA}_f \ZA \in \FT \ZA =\HO Z)$ varying smoothly with
$\ZA$ in $\JABG$. The definition of  $s_f \ZA$ is as follows.

   Choose an ordering $z_1,\ldots,z_d $ of points of the configuration $Z$ and let $\delta_i$ be the delta-function on $Z$ supported at $z_i$, for
$i=1,\ldots,d$. Decompose each $\delta_i$ according to the direct sum decomposition of
$\FT \otimes \OO_{\JABG}$ in (\ref{ordFT}) and let $\delta^0_i $ be its component
in $\HH^0 \ZA$. We now define
\begin{equation}\label{sfZA}
s^{\GA}_f \ZA =f(\delta^0_1,\ldots,\delta^0_d) \in \HO Z)
\end{equation}
by substituting $\delta^0_i$ for the indeterminate $T_i$, for $i=1,\ldots,d$. Since $f$ is a symmetric polynomial, $s^{\GA}_f \ZA$ is independent of the ordering chosen
and hence, defines a section of $\FT \otimes \OO_{\JABG}$ as $\ZA$ varies in $\JABG$. 
\\
\\
\indent
The perverse sheaves obtained by the two constructions above are similar to the ones in Theorem \ref{ps}. Namely, those are the intersection cohomology complexes
$IC(U_{\GA}, {\cal L})$
associated to $(U_{\GA}, {\cal L})$, where $U_{\GA}$ is a suitable Zariski open subset of $\JABG$ and ${\cal L}$ is the local system on $U_{\GA}$, modeled on the cohomology
ring of the Springer fibres, transplanted to $U_{\GA}$ using the maps in (\ref{D+}). In particular, these perverse sheaves are of representation theoretic nature.
Moreover, the Springer fibres above correspond to partitions of $d^{\prime}_{\GA}$ (see Theorem \ref{irrSn} for notation) and the shape of these partitions  is related, via
the representation theoretic version of Petri's method discussed in \S10.2, to equations
of configurations of $X$ parametrized by $\GA$.   Thus the perverse sheaves $IC(U_{\GA}, {\cal L})$
 can be viewed as a categorical way of writing equations of $\ZD_{\GA}$, the part of the universal scheme $\ZD$ lying over $\GA$, for every simple admissible component
$\GA$ in $\CS$. It should be clear now that the Jacobian $\JA$ together with the dictionary (\ref{ff-dict}) helps to put on a more solid footing our dream of recovering correspondences from
representation theoretic invariants.
 
\section{Concluding remarks and speculations}
The results discussed in {\bf Part II} give a convincing evidence
 that the Lie algebraic aspects of our Jacobian are useful in addressing various issues related to algebro-geometric properties of
configurations of points on surfaces. It also enables us to attach to the degree of the second Chern class of vector bundles such objects as irreducible
representations of symmetric groups and perverse sheaves. In fact we believe that the tools developed in \cite{[R2]} allow one to transfer virtually any object/invariant
of the geometric representation theory to the realm of smooth projective surfaces. Thus, for example, one should be able to have a version of Theorem \ref{irrSn}, where
the representations of the symmetric groups are replaced by the representations of the corresponding Hecke algebras as well as Affine Hecke algebras.

To our mind all these invariants fit into a sort of `secondary' type invariants for vector bundles in the sense of Bott and Chern 
in \cite{[B-C]}.
 Indeed, our construction begins by replacing the second Chern class of a bundle 
$\SE$ (of rank 2)
 by its geometric
realization, i.e. the zero-locus $Z$ of a suitable global section $e$ of $\SE$. This is followed by a distinguished orthogonal decomposition (\ref{ord-ZA-II}) of
the space of functions $\HO{Z})$ on $Z$. The decomposition gives rise to the Lie subalgebra ${\bf \tilde{\goth g}}(\SEE)$ of ${\bf \goth gl}(\HO{Z}))$ which is intrinsically associated
 to the pair\footnote{we return to the alternative, equivalent, notation of points in $\JA$ discussed in the last paragraph of \S2.1, i.e. $(\SEE) =\ZA$; in particular, with this matching
 the Lie algebra ${\bf \tilde{\goth g}}(\SEE) ={\bf \tilde{\goth g}}\ZA$.} 
$(\SEE)$. This Lie subalgebra could be viewed as the `secondary' structure Lie algebra associated to $\SE$. While the structure group 
(${\bf GL_2 (C)}$) with
its Lie algebra provide the topological invariants of $\SE$, i.e. its Chern classes, the secondary structure Lie algebra detects various algebro-geometric properties of the subscheme $Z$.
Thus Theorem \ref{tcd}, for example, can be interpreted as a statement of reduction of the secondary structure Lie algebra to a proper Lie subalgebra of
${\bf \goth gl}(\HO{Z}))$ (see (\ref{Lid})).  A geometric significance of such a reduction is the decomposition of $Z$ in (\ref{cd}). 
Furthermore, while the structure group and its Lie algebra
yield the Chern invariants of $\SE$ by evaluating the basic structure group-invariant polynomials on a curvature form of $\SE$, it is plausible to expect that our
 secondary Lie algebra
should provide many more representation theoretic invariants of $(\SEE)$ which would reflect  properties of geometric representatives of the Chern invariants of $\SE$.
The results presented in {\bf Part II} could be viewed as a confirmation of this heuristic reasoning.

 One can review the discussion in \S12.2 from the same vantage point. Namely, the procedure of sending a symmetric polynomial $f$ 
 to the section $s^{\GA}_f$ of $\FT \otimes \OO_{\JABG}$, described there,  is reminiscent of the construction of Chern classes of a vector bundle on a complex manifold, where one evaluates
(the elementary) symmetric polynomials on a curvature form to obtain the Chern forms representing the Chern classes in the cohomology ring of a manifold. 
In our situation the resulting value $s^{\GA}_f$ is a section of the sheaf of rings $\FT \otimes \OO_{\JABG}$. The direct sum decomposition (\ref{ordFT}) gives the decomposition 
$$
s^{\GA}_f =\sum^{\LG}_{p=0} s^{\GA,p}_{f}
$$
where $ s^{\GA,p}_{f}$ is the projection of $s^{\GA}_f $ onto the summand $\HH^p$ of the decomposition (\ref{ordFT}). In particular, $ s^{\GA,p}_{f}$ is a global section of $\HH^p$ and the collection 
$\{ s^{\GA,p}_{f}\}_{p=0,\ldots,\LG}$ can be regarded as the secondary Chern classes of vector bundles parametrized by $\JABG$.

Furthermore, using the dictionary (\ref{ff-dict}) the section $s^{\GA}_f$ can be elevated to a perverse sheaf either on $\JA$ or on the Hilbert scheme $\XD$. Thus one obtains a map
from the ring of symmetric functions to the category of perverse sheaves on the Hilbert schemes of $X$ and on its non-abelian Jacobians - a categorical version
of Chern-Weil homomorphism.


In the present work we pursued the direction of using the representation theory of reductive Lie algebras to study the geometry of surfaces. One could try to reverse the logic
and ask if linking representation theoretic objects to the geometry of the Hilbert schemes of surfaces could provide new insights or even new results in the representation theory.
A model for this line of inquiry is the work of M. Haiman, \cite{[Hai]}. In particular, one can ask\footnote{the question was asked by D. Kazhdan.}
 if Haiman's results could be recovered by applying our theory to $\PP^2$. At this stage we do not know the answer but the results in \cite{[R2]}, \S10.5, indicate that this is quite
plausible.

The results of Theorem \ref{ps} and Theorem \ref{LD} could be an indication that our constructions might be useful for questions of the categorification of representations
of simple Lie algebras of type ${\bf \goth sl}_n$.

The results of \S \ref{AffLie} and the discussion in \S \ref{Langlands} indicate a relation of our Jacobian to the Langlands duality. 
As suggested in \S \ref{Langlands}, a
formulation of the geometric Langlands Program for higher dimensional varieties could involve correspondences in the middle dimension. 
What we have in mind here is that
correspondences in the middle dimension could be taken as a geometric substitute for the Galois side of the Langlands duality.
Now
the very idea of the Jacobian as a tool to study correspondences goes back to A. Weil, \cite{[W]}. 
The scheme $\JA$ exhibits many promising features for being such a tool in the case of projective surfaces. 
 This together with the discussion in \S \ref{Langlands} can now be summarized as
the following triangular relation
\begin{equation}\label{tri}
\xymatrix{
&\JA \ar[dl] \ar[dr]&\\
*\txt{Correspondences\\ of X} \ar[rr] & &*\txt{Langlands Correspondence}} 
\end{equation}
A more detailed discussion of these interrelations will appear elsewhere but we hope that the results and ideas discussed in this survey will convince the reader
that the nonabelian Jacobian $\JA$ exhibits strong ties with the base of the above triangle.

 \begin{flushright}
 Igor Reider\\      
L'Universit\'e d'Angers
\\
reider@univ-angers.fr
\end{flushright}
\end{document}